\documentclass{amsart}
\usepackage {amssymb}
\usepackage {amsmath}
\usepackage{amsthm}
\usepackage{graphicx}
\usepackage {amscd}
\usepackage {epic}
\usepackage[alphabetic]{amsrefs}
\usepackage[colorlinks, citecolor=blue]{hyperref}

\newcommand{\spec}[0]{\operatorname{Spec}} 
\newcommand{\red}[0]{\operatorname{red}} 
\newcommand{\supp}[0]{\operatorname{Supp}}

\newcommand{\Tr}[0]{\operatorname{Tr}} 
\newcommand{\Pic}[0]{\operatorname{Pic}} 

\newcommand{\MO}{\mathcal{O}}

\newcommand{\Supp}{{\rm Supp}}
\newcommand{\Hom}{{\rm Hom}}
\newcommand{\sm}{{\rm sm}}

\newtheorem{thm}{Theorem}[section]

\newtheorem{lem}[thm]{Lemma}
\newtheorem{cor}[thm]{Corollary}

\newtheorem{prop}[thm]{Proposition}

\theoremstyle{definition}
\newtheorem{defn}[thm]{Definition}

\newtheorem{say}[thm]{}

\newtheorem{rem}[thm]{Remark}          

\newtheorem{notation}[thm]{Notation}   
  
\newtheorem{defn-thm}[thm]{Definition--Theorem}  
\newtheorem{defn-prop}[thm]{Definition--Proposition}  
\newtheorem{defn-lem}[thm]{Definition--Lemma}  

\theoremstyle{remark}

\input{diagrams.sty}

\begin{document}
\title{On base point freeness in positive characteristic}

\renewcommand{\subjclassname}{%
\textup{2010} Mathematics Subject Classification}
\subjclass[2010]{14E30, 13A35.}
\keywords{birational geometry, positive characteristic}

\author{Paolo Cascini}
\address{Department of Mathematics\\
Imperial College London\\
180 Queen's Gate\\
London SW7 2AZ, UK}
\email{p.cascini@imperial.ac.uk}

\author{Hiromu Tanaka}
\address{Department of Mathematics, Graduate School of Science, Kobe University, Kobe, 657-8501, Japan}
\email{tanakahi@math.kobe-u.ac.jp}

\author{Chenyang Xu}
\address{Beijing International Center of Mathematics Research, 5 Yiheyuan Road, Haidian District, Beijing, 100871, China}
\email{cyxu@math.pku.edu.cn}

\begin{abstract}
We prove that if 
 $(X,A+B)$ is a pair defined over an algebraically closed field of positive characteristic  such that   $(X,B)$ is strongly $F$-regular, $A$ is  ample and $K_X+A+B$ is strictly nef, then $K_X+A+B$ is ample. Similarly, we prove that for a log pair $(X,A+B)$ with $A$ being ample and $B$ effective, $K_X+A+B$ is big if it is nef and of maximal nef dimension. As an application, we establish a rationality theorem for the nef threshold and various results towards the minimal model program in dimension three in positive characteristic.
\end{abstract}

\date{\today}
\maketitle{}

\tableofcontents

\section{Introduction}

One of the main objectives of the minimal model program is the study of the linear system associated to an adjoint divisor. For example,  in characteristic 0, we 
have a good understanding of the linear system given by a multiple of a  $\mathbb Q$-divisor $L$ which is  the sum of the canonical divisor and an ample $\mathbb Q$-divisor (e.g. see \cite{KMM85}, \cite{KM98}, \cite{bchm} and the 
references therein). A fundamental tool in birational geometry is Kawamata's base point free theorem which asserts that if such a $\mathbb Q$-divisor $L$ is nef then it is semiample (see \cite{KM98}). 

Because of the failure of the Kodaira vanishing theorem in positive characteristic, 
Kawamata's base point free theorem and its generalisations are not known to hold in this case. The aim of this paper is to present a new approach  to the base point free theorem in positive characteristic. 
We prove an special case of this result as well as several results which, in characteristic $0$, are known to follow from the base point free theorem.


\subsection{Strictly nef divisors} 
We first study strictly nef adjoint divisors, with possibly real coefficients. Recall that  a $\mathbb R$-Cartier  $\mathbb R$-divisor $L$ on a  proper variety $X$ is said to be {\it strictly nef} if its intersection with any curve on $X$ is positive. 
Mumford has constructed the first example of a strictly nef divisor which is not ample (see \cite[Example 10]{Hartshorne70}).  
See  \cite[Remark 3.2]{MS95} for a similar example in positive characteristic.
However, we show:
\begin{thm}\label{t-main}
Let $(X,B)$ be a strongly $F$-regular pair defined over an algebraically closed field $k$ of characteristic $p>0$, where $B$ is an effective $\mathbb R$-divisor. Assume that  $A$ is an ample $\mathbb R$-divisor such that  
 $K_X+A+B$ is strictly nef. 
 Then $K_X+A+B$ is ample.
\end{thm}

\medskip 
From Theorem \ref{t-main}, we immediately obtain the following result:
\begin{cor}\label{c-main}
 Let $(X,\Delta)$ be a strongly $F$-regular projective pair 
with an effective $\mathbb R$-divisor $\Delta$ over an algebraically closed field $k$ of characteristic $p>0$ such that $K_X+\Delta$ is big and strictly nef. 
 Then $K_X+\Delta$ is ample. 
\end{cor}

In addition, we  obtain the following result on the rationality of the nef threshold: 
\begin{thm}\label{t-rational}
Let $(X,B)$ be a strongly $F$-regular pair defined over an algebraically closed field of characteristic $p>0$, where $B$ is an effective $\mathbb Q$-divisor. Assume that 
$K_X+B$ is not nef and $A$ is an ample $\mathbb Q$-divisor. Let 
$$\lambda:=\min\{t> 0\mid K_X+B+tA\text{ is nef }\}.$$

Then there exists a curve $C$ in $X$ such that $(K_X+\lambda A+B)\cdot C=0$. In particular, $\lambda$ is a rational number. 
\end{thm}

When $X$ is smooth and $B=0$, the results follow from Mori's cone theorem \cite{Mori82}. We note that the assumption that $(X,B)$ is strongly $F$-regular is analogous  but more restrictive than the assumption that $(X,B)$ is klt.  In fact, in characteristic 0, all these statements are direct consequences of Kawamata's base point free theorem as we know that if $(X,B)$ is a projective klt pair such that $B$ is big and $K_X+B$ is nef then $K_X+B$ is indeed semi-ample (see e.g. \cite[3.3]{KM98}).
 
 In positive characteristic, since \cite{HH90} new techniques involving the Frobenius map have been developed to establish 
 the positive characteristic analogs of many of the results, which in characteristic $0$, are traditionally deduced from vanishing theorems. 
Very roughly, this is the general strategy that we follow in this paper as well.

On the other hand, the techniques used to  prove the above results were inspired by an earlier attempt of the second author to prove Fujita's conjecture, which in turn was inspired by the proof of the effective base point free theorem in characteristic zero, by Angehrn and Siu \cite{AS95}. 
In their paper, the authors construct  zero-dimensional subschemes which are minimal log canonical centres for a suitable pair and using Nadel's vanishing theorems they are able to extend non-trivial sections to the whole variety. In positive characteristic, using the idea of twisting by Frobenius, the analogue would be to construct zero dimensional $F$-pure centres and use $F$-adjunction (see \cite{Sch09} for more details).  Unfortunately there is a technical issue due to the index of the adjoint divisor, which we are not able to deal with, in general. Therefore, instead of using one divisor to cut the centre, we study the trace map for all the powers of Frobenius and assign different coefficients for each of these. 
For this reason, we introduce the use of $F$-{\em threshold functions} to replace the classical $F$-pure threshold and obtain a zero dimensional subscheme from which we can lift sections (see Subsection \ref{s-cut} and \ref{s-induction} for more details).    Theorem \ref{t-main} and Theorem \ref{t-rational} are proven in Section \ref{s-extend}.

\medskip 

\subsection{Divisors of maximal nef dimension}  Using the same methods as above but cutting at two very general points, we study adjoint divisors of maximal nef dimension.  
More specifically, given a log pair $(X,B)$ such that $K_X+B$ is nef, the nef reduction map associated to $K_X+B$ (see Subsection \ref{ss-nrm} for the definition) has proven to be a powerful tool to approach the Abundance conjecture  (e.g. see \cite[Section 9]{BDPP04} and \cite{a04} for more details). 
 Recall that a divisor over a proper variety $X$ is said to be of {\em  maximal nef dimension} if its intersection with any movable curve in $X$ is positive (see Subsection \ref{s-nc} for the definition of movable curve). Thus, we obtain the following weak version of the base point free theorem: 

\begin{thm}\label{t-main2}
Let $X$ be a normal projective variety over an algebraically closed field  of characteristic $p>0$. Assume that  $A$ is an ample $\mathbb R$-divisor and $B\ge 0$ is a $\mathbb R$-divisor such that 
$K_X+B$ is $\mathbb R$-Cartier and 
 $K_X+A+B$ is nef and of maximal nef dimension.  Then $K_X+A+B$ is big.
\end{thm}
Note that the previous theorem does not require any assumption on the singularities of the pair $(X,B)$, nor on the coefficients of $B$.

As an application, we obtain the following result on the extremal ray associated to a nef but not big adjoint divisor:

\begin{cor}\label{c-rational-cover}
Let $X$ be a normal projective variety, defined over an algebraically closed field of characteristic $p>0$. 
Assume that $A$ is an ample $\mathbb R$-divisor, $B\ge 0$ is an $\mathbb R$-divisor such that $K_X+B$ is $\mathbb R$-Cartier and $L=K_X+A+B$ is nef and not big. Assume that
$$\overline{NE}(X)\cap L^{\perp}=R$$
is an extremal ray of $\overline{NE}(X)$. 

Then $X$ is covered by rational curves $C$ such that $[C]\in R$ and $$-(K_X+B)\cdot C\le 2\dim X.$$
\end{cor}

Theorem \ref{t-main2} and Corollary \ref{c-rational-cover} are proven in Section \ref{s-max}.

\begin{rem}We were informed by J. M$^{\rm c}$Kernan that Theorem \ref{t-main2} and Corollary \ref{c-rational-cover} were  independently obtained  by him using  different methods \cite{McK}.
\end{rem}

\medskip

\subsection{Threefolds}
We now  focus on the study of three dimensional projective varieties. 
We first prove the following version of the cone theorem:
\begin{thm}\label{t-cone}
Let $X$ be a $\mathbb Q$-factorial projective  threefold defined over an algebraically closed field of
characteristic  $p>0$. Let $B$ be an effective $\mathbb Q$-divisor on $X$ whose coefficients are strictly less than one. 
Assume that $K_X+B$ is not nef. 
Then there exist an ample $\mathbb Q$-divisor  $A$ such that  $K_X+A+B$ is not nef and finitely many curves  $C_1, \cdots , C_r$ on $X$ such that 

$$\overline{NE}(X)=\overline{NE}(X)_{K_X+A+B \ge 0}+\sum_{i=1}^r \mathbb R_{\geq 0}[C_i].$$

\end{thm}

By combining our results with previous ones \cite{Kol91,Keel99,HX13}, we obtain a weak version of the minimal model program for three dimensional varieties:
\begin{thm}\label{t-3mmp}
Let $X$ be a $\mathbb{Q}$-factorial terminal projective threefold defined over an algebraically closed field of characteristic $p>5$. Then there exists a  $K_X$-negative birational contraction $f:X\dasharrow Y$ to a  $\mathbb{Q}$-factorial terminal projective threefold such that one of the following is true:
\begin{enumerate}
\item if $K_X$ is pseudo-effective, then $K_Y$ is nef;
\item if $K_X$ is not pseudo-effective, then there exist a $K_Y$-negative extremal ray $R$ of $\overline{NE}(Y)$ and a surjective morphism $g: Y\to Z$ to a normal projective variety $Z$ such that 
$\dim Y>\dim Z$, 
$g_*\mathcal O_Y=\mathcal O_Z$ and for every curve $C$ in $Y$, $g(C)$ is a point if and only if $[C]\in R$.
\end{enumerate}
\end{thm}

Theorem \ref{t-cone} and Theorem \ref{t-3mmp} are proven in Subsection \ref{s-cone}.

\medskip 

We also prove the following version of the base point free theorem in  dimension three, under some assumptions on the coefficients of the boundary:

\begin{thm}\label{t-abundance}
Let $(X,B)$ be a projective three dimensional log canonical pair defined over an algebraically closed field of characteristic $p>0$, for some big $\mathbb Q$-divisor $B\ge 0$ such that  $K_X+B$ is nef.  Assume that $p>\frac 2 a$ for any coefficient $a$ of $B$.
\begin{enumerate}
\item If $K_X+B$ is not numerically trivial, then 
$$\kappa(X,K_X+B)=\nu(X,K_X+B)=n(X,K_X+B).$$
\item If $\kappa(X,K_X+B)=1$ or $2$, then $K_X+B$ is semiample.
\item If $k=\overline{\mathbb{F}}_p$, and all coefficients of $B$ are strictly less than 1, then $K_X+B$ is semiample.
\end{enumerate}
\end{thm}

Note that if $(X,B)$ is a three dimensional projective log pair such that $K_X+B$ is big and nef then Keel proved a version of the base point free theorem which allows the target space to be an algebraic space (cf. \cite[Theorem 0.5]{Keel99}).

Besides using Theorem \ref{t-main2}, the main tool used to prove Theorem \ref{t-abundance} is a canonical bundle formula for fibrations of relative dimension one. The proof is contained in Subsection \ref{s-bpf}.

\medskip 

\noindent \textbf{Acknowledgement.} 
The first author is partially supported by an EPSRC Grant. The third author is partially supported by the grant `The Recruitment Program of Global Experts'. Part of this work was completed while the first author was visiting the Beijing International Center of Mathematics Research. He would like to thank the third author for his very generous hospitality. 

We would like to thank O. Fujino, Y. Gongyo, C. Hacon, J. Koll\'ar, A.F. Lopez, J. M\textsuperscript cKernan, Z. Patakfalvi and K. Schwede for several very useful discussions. We would also like to thank the referee for carefully reading the paper and for many helful suggestions. 

\section{Preliminary results}

\subsection{Notation and conventions}\label{s-nc} We work over an algebraically closed field $k$ of positive characteristic $p$. If $K$ is a field, we denote by $\overline K$ its algebraic closure. By abuse of notation, we will often write $K$ instead of ${\rm Spec } K$. 

A {\em variety} $X$ is an integral scheme which is separated and of finite type over $k$. 
A {\em curve} is a one dimensional variety.  A curve $C$ in a variety $X$ is said to be {\em movable} if it is a member of an algebraic family $\mathcal{C}/T=(C_t)_{t\in T}$ parametrized by a variety $T$ and such that $\mathcal{C}\to X$ is dominant. 

When the ground field is uncountable, by {\em a very general point $x\in X$}, we mean a point $x$ which is in a subset $U$ given by the complement  of a countable union of proper subvarieties. By a pair of very general points, we mean $(x,y)\in U\times U$.

Let $K\in \{\mathbb Q, \mathbb R\}$. A $K$-{\em line bundle} $L$ on a proper scheme $X$ 
is an element of the group $\Pic(X)\otimes K$. We will use the additive notation on this group. A $K$-line 
bundle $L$ on $X$ is said to be {\em nef} (respectively {\em strictly nef}, {\em numerically trivial}) if 
$L\cdot C\ge 0$ (respectively $>0$, $=0$) for all the curves $C$ in $X$. The $K$-line bundle $L$ is said to 
be of {\em maximal nef dimension} if $L\cdot C>0$ for all the movable curves $C$ in $X$.

If $X$ is a normal variety, we denote by ${\rm Div}_{\mathbb R}(X)$ the vector space of $\mathbb R$-Cartier $\mathbb R$-divisors of $X$, by $N_1(X)$ the vector space of $1$-cycles on $X$ up to numerical equivalence, and by
$\overline{NE}(X)\subseteq N_1(X)$ the closure of the convex cone generated by the classes of effective 
$1$-cycles in $X$. If $L$ is an $\mathbb R$-Cartier $\mathbb R$-divisor on $X$, we denote by $L^\perp\subseteq N_1(X)$ the 
set of $1$-cycles $C$ on $X$  such that $L\cdot C=0$ and by $\overline{NE}(X)_{L\ge 0}$ the set of $1$-cycles $C\in 
\overline{NE}(X)$ such that $L\cdot C\ge 0$. 
 Given any 
$\mathbb{R}$-Cartier $\mathbb{R}$-divisor $D$ and an ample $\mathbb{R}$-divisor $H$, we define {\it the nef threshold}  of $D$ 
with respect to $H$ to be $$\lambda =\min \{t\ge 0\mid D+tH \text{ is nef} \}.$$

Given a $\mathbb Q$-line bundle $L$ on a proper variety $X$, we denote by $\kappa(X,L)$ its Iitaka dimension and we define its {\em volume} as 
$${\rm vol}(X,L)=\limsup_{m\to \infty} \frac{n!~h^0(X,mL)}{m^n}
$$
where $n$ is the dimension of $X$ and $m$ is taken to be sufficiently divisible (see \cite[\S 2.2.C]{Laz} for more details).  
If $L$ is a nef $\mathbb{R}$-line bundle, we denote by $\nu(X,L)$ the numerical dimension of $L$, i.e.  
$$\nu(X,L)=\max\{m\ge 0~|~ L^m \not\equiv 0\}$$ 
where we denote by $L^m$ the $m$-th self intersection of $L$. 

We refer to \cite{KM98} for the classical definitions of singularities (e.g., {\it klt, log canonical}) appearing in the minimal model program, except for the fact that in our definitions we require the pairs to have {\it effective} boundaries.    In addition, given a $\mathbb Q$-divisor $B$ on a normal variety $X$  such that $K_X+B$ is $\mathbb Q$-Cartier, we say that the pair $(X,B)$ is {\em sub log canonical} if $a(E,X,B)\ge -1$ for any geometric valuation $E$ over $X$. 

Given a variety $X$, we denote by $F\colon X\to X$ the  absolute Frobenius morphism. We refer to Definition \ref{char-F-sing} for the definition of a {\em strongly $F$-regular pair} and a {\em sharply $F$-pure} pair. 
If $Z$ is a closed subscheme of a projective variety $X$, then the scheme-theoretic inverse image
$$Z^{[e]}:=(F^e)^{-1}(Z)$$
is a closed subscheme of $X$ defined by the ideal $I^{[p^e]}_Z$, so that if $I_Z$ is locally generated by $f_1,\dots,f_k$ then $I^{[p^e]}_Z$ is locally defined by $f_1^{p^e},\dots,f_k^{p^e}$.

\subsection{Preliminaries} We begin with the following well known results. 

\begin{lem}\label{s-nef}
A nef $\mathbb R$-line bundle $L$ on a projective variety $X$ is strictly nef if and only if $L|_{V}$ is not numerically trivial for any subvariety $V\subseteq X$ with $\dim(V)\ge 1$.
\end{lem}
\begin{proof}Pick an ample divisor $H$ on $X$. Assume that $L$ is strictly nef and that $V\subseteq X$ is a subvariety. Then 
$$L|_V\cdot (H^{\dim V-1}\cdot V)>0.$$ Thus, $L|_V$ is not numerically trivial. The converse is trivial.   
\end{proof}

\begin{lem}\label{l-nef2}
Assume that $X$ is a projective variety defined over an uncountable algebraically closed field. Let $L$ be an  $\mathbb R$-line bundle of maximal nef dimension on $X$.

Then, for  a very general point $x\in X$,  $L|_{V}$ is not numerically trivial for any subvariety $V\subseteq X$ such that 
$\dim V\ge 1$ and $x\in V$. 
\end{lem}
\begin{proof}Cutting by hyperplanes, it suffices to prove that for a very general point $x$ and for any irreducible curve $C$ through $x$, the restriction $L|_C$ is not numerically trivial. 

Let ${\rm Univ}_1\to {\rm Chow}_1$  be the universal family over the Chow variety parameterizing 1-dimensional cycles.  Note that the set of  non-movable curves $C\subseteq X$ is parametrized by a countable union of subvarieties $W\subseteq {\rm Chow}_1$ such that ${\rm Univ}_W\to X$ is not dominant. Let $x$ be a very general point which is not contained in the union of the closures of the image  of each component of ${\rm Univ}_W$.

Then, the lemma  follows from the fact that $L$ is of maximal nef dimension and any curve $C$ through $x$ is a movable curve.    
\end{proof}

We need the following ampleness criterion in Section \ref{s-extend}:

\begin{lem}\label{l-nef3}
Let $L$ be a  strictly nef $\mathbb{R}$-Cartier $\mathbb R$-divisor on a normal projective variety $X$. Assume that 
for every closed point $x \in X$, we may write $L \sim_{\mathbb{R}}L_x$ where 
$L_x$ is an effective $\mathbb R$-divisor whose support does not contain $x$. Then $L$ is ample.
\end{lem}
\begin{proof} By the Nakai-Moishezon theorem (for $\mathbb{R}$-divisors, see \cite{CP90}), we only need to check for any subvariety $Z$ of $X$, $L^{\dim Z}\cdot Z>0$. 
By induction on the dimension, we can assume that for any proper subvariety $Y\subsetneq X$, if $\nu\colon \overline Y\to Y$ is the normalisation of $Y$, then  $L^{\dim Y}\cdot Y=(\nu^*L|_Y)^{\dim Y}>0$.

By  assumption, we can write $L=\sum^q_{i=1} c_i L_i$  for some positive numbers $c_1,\dots,c_q$ and distinct prime divisors $L_1,\dots,L_q$. Therefore,  
$$L^n=\sum_{i=1}^q c_i (L|_{L_i})^{n-1}>0,$$
and the claim follows.  
\end{proof}

Similarly, we need the following bigness criterion  in Section \ref{s-max}: 

\begin{lem}\label{l-nef4}
Let $X$ be a normal projective variety, defined over 
an uncountable algebraically closed field. 
Let $L$ be a nef $\mathbb R$-Cartier $\mathbb R$-divisor. 
Assume that, for any very general points $x, y\in X$,  there exists an effective $\mathbb R$-Cartier $\mathbb R$-divisor $L_{x,y}\sim_{\mathbb R}L$
 such that $x\in \Supp ~L_{x,y}$ and  $y\not\in \Supp~ L_{x,y}$. 
Then $L$ is big. 
\end{lem}

\begin{proof}
It is sufficient to show that $L^n>0$. 
Fix a very general point $x\in X$. Then  
we can find an effective $\mathbb R$-Cartier $\mathbb R$-divisor $L_1 \sim_{\mathbb{R}}L$ containing $x$ in its support. We may write $L_1=fF+G$ where $F$ is a prime divisor such that $x\in F$, $f$ is a positive number and $G$ is an effective $\mathbb R$-divisor which does not contain $F$ in its support. 
Note that, since $x\in F$,  if $\nu\colon \overline F\to F$ is the normalisation, then $\nu^*(L|_{F})$ satisfies the same properties as $L$. 
Thus, by induction on the dimension, we obtain  
$$L^n\ge L^{n-1}\cdot fF=f(L|_{F})^{n-1}=f(\nu^*L|_F)^{n-1}>0.$$ 
Thus, the claim follows. 
\end{proof}

\medskip

\subsection{The trace map of Frobenius}\label{ss-trace}

All the results in this section are essentially contained in the fundamental work \cite{Sch09}.  We include them for the reader's convenience. 

\begin{defn-prop}\label{trace-def}
Let $X$ be a normal variety, 
let $D$ be an effective divisor on $X$ and let $e$ be a positive integer. 
Then we can define an $\mathcal O_X$-module homomorphism 
\begin{eqnarray*}
{\rm Tr}_X^e(D)\colon F_*^e(\mathcal O_X(-(p^e-1)K_X-D))\to \mathcal O_X
\end{eqnarray*}
which satisfies the following commutative diagram of $\mathcal O_X$-modules 
$$\begin{CD}
F_*^e(\mathcal O_X(-(p^e-1)K_X-D))@>{\rm Tr}_X^e(D) >> \mathcal O_X\\
@V\theta V\simeq V @VV\simeq V\\
\mathcal Hom_{\mathcal O_X}(F_*^e(\mathcal O_X(D)), \mathcal O_X)
@>(F^e(D))^*>>\mathcal Hom_{\mathcal O_X}(\mathcal O_X, \mathcal O_X).
\end{CD}$$
\end{defn-prop}

The result above has appeared in the literature before (e.g. see \cite[Section 2]{Sch09}, \cite[Section 2]{Tan13}). We provide a proof here for the sake of completeness.

\begin{proof}
Consider the composition map 
$$F^e(D)\colon \mathcal O_X\overset{F^e}\to F_*^e\mathcal O_X\hookrightarrow F_*^e(\mathcal O_X(D)).$$
Apply the contravariant functor $\mathcal Hom_{\mathcal O_X}(-, \mathcal O_X)$: 
 $$(F^e(D))^*:\mathcal Hom_{\mathcal O_X}(F_*^e(\mathcal O_X(D)), \mathcal O_X)
\to \mathcal Hom_{\mathcal O_X}(F_*^e\mathcal O_X, \mathcal O_X)
\to \mathcal Hom_{\mathcal O_X}(\mathcal O_X, \mathcal O_X). $$
We want to show that there exists an $\mathcal O_X$-module isomorphism: 
{$$\mathcal Hom_{\mathcal O_X}(F_*^e(\mathcal O_X(D)), 
\mathcal O_X)\simeq F_*^e(\mathcal O_X(-(p^e-1)K_X-D)).$$} 
Note that both  coherent sheaves are reflexive. Denote by $i\colon X^{\sm}\hookrightarrow X$ the open embedding of the smooth locus of $X$.  
We have 
{$$\mathcal Hom_{\mathcal O_X}(F_*^e(\mathcal O_X(D)),\mathcal{O}_X)\cong i_*\mathcal Hom_{\mathcal O_{X^{\sm}}}(F_*^e\mathcal O_{X^{\sm}}(D|_{X^{\sm}}),\mathcal{O}_{X^{\sm}})$$} 
and
$$F_*^e(\mathcal O_X(-(p^e-1)K_X-D))\cong i_* F_*^e(\mathcal O_{X^{\sm}}(-(p^e-1)K_{X^{\sm}}-D|_{X^{\sm}})).$$
Therefore, replacing $X$ by its smooth locus, we may assume that $X$ is smooth. 
By the duality theorem for finite morphisms, 
we obtain the following $\mathcal O_X$-module isomorphism 
\begin{eqnarray*}
&\theta:&\mathcal Hom_{\mathcal O_X}(F_*^e(\mathcal O_X(D)), \mathcal O_X)\\
&\simeq& 
\mathcal Hom_{\mathcal O_X}(F_*^e(\mathcal O_X(D)), \omega_X)\otimes (\omega_X)^{-1}\\
&\simeq& 
F_*^e\mathcal Hom_{\mathcal O_X}(\mathcal O_X(D), \omega_X)\otimes (\omega_X)^{-1}\\
&\simeq&F_*^e\mathcal Hom_{\mathcal O_X}(\mathcal O_X(D), \mathcal O_X((1-p^e)K_X))\\
&\simeq&F_*^e(\mathcal O_X(-(p^e-1)K_X-D)).
\end{eqnarray*}
Thus, the claim follows. 
\end{proof}

\begin{prop}\label{trace-vs-split}
Let $X$ be a normal variety, 
 and let $D$ be an effective divisor on $X$.
Fix a  positive integer $e$ and a scheme-theoretic point $x\in X$.
Then, the following assertions are equivalent:
\begin{enumerate}
\item ${\rm Tr}_X^e(D)$ is surjective at $x$;
\item The $\mathcal O_{X, x}$-module homomorphism 
$$(F^e(D))_x\colon \mathcal O_{X, x}\overset{F^e}\to F^e_*\mathcal O_{X, x}\hookrightarrow F^e_*(\mathcal O_{X,x}(D))$$
splits. 
\end{enumerate}
\end{prop}

\begin{proof}
Assume (1). 
Then, there exists $\varphi \in \Hom_{\mathcal O_{X, x}}(F_*^e(\mathcal O_{X, x}(D)), \mathcal O_{X,x})$
such that 
$$(F^e(D))^*(\varphi)={\rm id}_{\mathcal O_{X,x}}$$ 
Thus, $\varphi$ gives the required splitting. 

Assume (2). 
Then, there exists 
$\varphi \in \Hom_{\mathcal O_{X, x}}(F_*^e(\mathcal O_{X, x}(D)), \mathcal O_{X,x})$
such that 
$(F^e(D))^*(\varphi)={\rm id}_{\mathcal O_{X,x}}.$ 
This implies the required surjectivity. 
\end{proof}

\begin{defn}\label{char-F-sing}
Let $X$ be a normal variety and let $B$ be an effective $\mathbb R$-divisor 
such that $K_X+B$ is $\mathbb R$-Cartier. 
Fix a closed point $x\in X$. 
\begin{enumerate}
\item{A pair $(X, B)$ is {\em strongly $F$-regular} at $x$ if, for every effective divisor $E$, 
there exists a positive integer $e$ such that 
${\rm Tr}_X^e(\ulcorner(p^e-1)B\urcorner+E)$ is surjective at $x$.}
\item{A pair $(X, B)$ is {\em sharply $F$-pure} at $x$ if  
there exists a positive integer $e$ such that 
${\rm Tr}_X^e(\ulcorner(p^e-1)B\urcorner)$ is surjective at $x$.}
\end{enumerate}
\end{defn}

\begin{rem}\label{F-sing-rem}
\begin{enumerate}
\item{By Proposition~\ref{trace-vs-split}, 
if $B$ is a $\mathbb Q$-divisor, then the above definition coincides with the one in \cite[Definition~2.7]{Sch09}. }
\item{If $(X,B)$ is strongly $F$-regular and $E$ is an effective divisor, then our definition implies that there exists a $\mathbb{Q}$-divisor $B'\ge B$, such that the map 
${\rm Tr}_X^e(\ulcorner(p^e-1)B'\urcorner+E)$ is surjective.
Thus for any  effective Cartier divisor $D$, we choose $E$ in Definition \ref{char-F-sing} to be a sufficiently large effective divisor whose supports  contains ${\rm Supp }(B'+D)\cup {\rm Sing}(X)$. Applying \cite[Theorem 3.9]{SS10}, we obtain that for a sufficiently small number $\varepsilon >0$, we have that $(X,B'+\varepsilon D)$ is strongly $F$-regular, which implies $(X,B+\varepsilon D)$ is also strongly $F$-regular as well. 
}
\item By abuse of notation, we will often denote ${\rm Tr}_X^e(D)$ simply by ${\rm Tr}^e$. 
\end{enumerate}
\end{rem}

\medskip

\subsection{ Nef reduction map}\label{ss-nrm}

We now  recall the main result of \cite{BC02}, which allows us to study nef line bundles on a projective variety which are not of maximal nef dimension. 

\begin{thm}[Nef reduction map]\label{t_nrm}
Let $X$ be a normal  projective variety defined over an uncountable algebraically closed field $k$, and $L$ be a nef $\mathbb R$-line bundle. Then there exists an open set $ U\subseteq X$ and a proper morphism $\varphi\colon U\to V$, such that $L$ is numerically trivial on a very general fibre $F$ of $\varphi$ and for a very general point $x$, we have that $L\cdot C=0$ if and only if $C$ is contained in  the fibre of $\varphi$ containing $x$.
\end{thm}

\begin{proof} The theorem follows from the main result in \cite{BC02}. Although the result there is stated only for line bundles on complex projective varieties, the same proof works for $\mathbb R$-line bundles on any variety defined over an uncountable algebraically closed field. 
\end{proof}

It follows from the previous theorem that if $L$ is a nef line bundle on a normal projective variety $X$ defined over an algebraically closed field $k$, we can define the  {\it  nef dimension $n(X,L)$} as the dimension  of the variety $V$ in  Theorem \ref{t_nrm}, after first possibly applying a base change so that $X$ is defined over an uncountable field $K\supseteq k$. It is clear that this definition does not depend on the choice of $K$. Note that $L$ is of maximal nef dimension if and only if $n(X,L)=\dim X$ and in general we have the inequalities:
$$\kappa(X,L)\le \nu(X,L)\le n(X,L)\le \dim X.$$ 
(see \cite[Proposition 2.8]{BC02}).

\section{Creating isolated centres}
In this section, we aim to develop the method of creating isolated centres. As we mentioned, our approach is different from the standard one, because instead of studying one threshold, we track a sequence of thresholds. In Subsection \ref{s-cut}, we study how to cut out $(d-1)$-dimensional centers from $d$-dimensional centers. In Subsection \ref{s-induction}, we establish the induction process for all $d$.

\subsection{Cutting subschemes}\label{s-cut}
 
In this section, we always assume that  $X$ is a projective variety defined over an algebraically closed field of  characteristic $p>0$. Our goal is to construct zero-dimensional subschemes of $X$ from which we can lift sections. These methods were inspired by the proof of the effective base point free theorem in characteristic zero, by Angehrn and Siu \cite{AS95}.

The following result will allow us to create isolated $F$-pure centres. 

\begin{prop}\label{q-cut}
Fix $a \in \mathbb{ N}$. Let $X$ be a projective variety. Let $A$ be an ample $\mathbb R$-line bundle on $X$ and $L$ a nef $\mathbb R$-line bundle on $X$. Let $x \in X$ be a closed point and let $ W$ be a proper closed subscheme of $X$. Assume $\dim_x W \ge 1$. 

Then there exist a positive integer  $\lambda_0$ and an ample $\mathbb R$-line bundle $A'$ on $X$ such that 
\begin{enumerate}
\item{$A-A'$ is ample, }
\item{$\lambda_0L+A'$ is a $\mathbb Q$-line bundle, and }
\item{for any sufficiently divisible $l>0$, there exists 
$$ t \in H^0(X, l(\lambda_0L+A') \otimes (m^{al}_x +  I_W))$$
such that $t|_V\neq 0$ for every irreducible component $V$ of $W^{\red}$ 
such that $L|_V$ is not numerically trivial.}
\end{enumerate}
\end{prop}

Before we proceed with the proof of Proposition \ref{q-cut}, we first need some preliminary results. 

\begin{lem}\label{l-cut1}
Fix $a \in \mathbb{ N}$. Let $X$ be a projective variety and let $V$ be a closed subvariety of $X$. Let $A$ be an ample $\mathbb R$-line bundle on $X$ and $L$ a nef $\mathbb R$-line bundle on $X$ such that $L|_V$ is not numerically trivial. Let $x \in V$ be a closed point and let $ W'$ be a proper closed subscheme of $V$. 

Then there exist a positive integer  $\lambda_0$, such that for any integer  $\lambda\ge \lambda_0$ there exist an 
ample $\mathbb R$-line bundle $A_{\lambda}$ on $X$ and a positive integer $q_{\lambda}$ such that 
\begin{enumerate}
\item $A-A_{\lambda}$ is ample,
\item $q_{\lambda}(\lambda L+A_{\lambda})$ is a line bundle, and 
\item  for every positive integer $l$, we have 
$$H^0(V, lq_{\lambda}(\lambda L+A_{\lambda})\otimes (m_x^{alq_{\lambda}} \cap I_{W'} ))\neq 0. $$
\end{enumerate}
\end{lem}

\begin{proof}  
For every positive integer  $\lambda$, 
we can find an ample $\mathbb R$-line bundle $A_{\lambda}$ such that 
\begin{itemize}
\item $A-A_{\lambda}$ is ample,
\item $\lambda L+A_{\lambda}$ is a $\mathbb Q$-line bundle, and
\item $A_{\lambda}-\frac{1}{2}A$ is ample.
\end{itemize}
Indeed, we can find such an $\mathbb R$-line bundle $A_{\lambda}$ by perturbing $\frac{3}{4}A$. 
Let $r$ denote the dimension of $V$. Fix a closed point $x\in V$. 
Let 
$$H_x(l)={\rm length}_x(\mathcal{O}_{x,V}/m^{l+1}_x)\qquad (l\gg 0)$$ be the Hilbert-Samuel function (see \cite[Section 12.1]{Eisenbud95}).   
Then $H_x(l)$ is of the form
$$H_x(l)=\frac{e_x\cdot l^r}{r!}+(\mbox{lower terms}),$$
where $e_x$ is the multiplicity of $x\in V$.
We have
\begin{eqnarray*}
\lim_{\lambda\to \infty} {\rm vol}(V,\lambda L+A_{\lambda})
&=&\lim_{\lambda\to \infty} (\lambda L|_V+A_{\lambda}|_V)^r\\
&\geq&\lim_{\lambda\to \infty} \lambda rL|_V\cdot (A_{\lambda}|_V)^{r-1}\\
&\geq& \lim_{\lambda\to \infty} \lambda rL|_V\cdot (\frac{1}{2}A|_V)^{r-1}\\
&=& \infty.
\end{eqnarray*}
Hence we can find a sufficiently large integer $\lambda_0$ such that 
$${\rm vol}(V,\lambda L+A_{\lambda})> e_x\cdot {a^r}$$
for every $\lambda\ge \lambda_0.$ 

Therefore, if $\lambda\ge \lambda_0$ and $l$ is sufficiently divisible, we have
$$\begin{aligned}
& h^0(V,l(\lambda L+A_{\lambda}) \otimes (m_x^{al} \cap I_{W'} ))\\
&\ge  ~ h^0(V,l(\lambda L+A_{\lambda}))-h^0(V,\MO_V/(m_x^{al} \cap I_{W'})\otimes l(\lambda L+A_{\lambda}))\\ 
&\ge  ~ h^0(V,l(\lambda L+A_{\lambda}))-h^0(V,\MO_V/m_x^{al}  \otimes l(\lambda L+A_{\lambda}))-h^0(W',  l(\lambda L
+A_{\lambda}))\\
&= ~ \frac{l^r{\rm vol}(V,\lambda L+A_{\lambda})}{r!}-\frac{e_x\cdot (al)^r}{r!}+ (\mbox{lower terms})\\
&= ~ \frac{{\rm vol}(V,\lambda L+A_{\lambda})-e_x\cdot a^r}{r!}l^r+ (\mbox{lower terms})\\
&\to ~ \infty\,\,\,(\text { if }l\to \infty~).
\end{aligned}$$
Thus, the claim follows. 
\end{proof}

\begin{lem}\label{l-cut2} 
Fix $a \in \mathbb{ N}$. Let $X$ be a projective variety and let $W$ be a reduced closed subscheme of $X$. Let $A$ be an ample $\mathbb R$-line bundle on $X$ and $L$ a nef $\mathbb R$-line bundle on $X$. Let $x \in W$ be a closed point. If $\dim_x W \ge 1$, then there exists a positive integer  $\lambda_0$ and an ample $\mathbb R$-line bundle $A'$ on $X$ such that 
\begin{enumerate}
\item{$A-A'$ is ample, }
\item{$\lambda_0L+A'$ is a $\mathbb Q$-line bundle, and }
\item{for any sufficiently divisible $l>0$,  there exists 
$$s\in H^0(W,l(\lambda_0L+A')\otimes m_x^{al})$$
such that $s|_V\neq 0$ for every irreducible component $V$ of $W$ such that 
$L|_V$ is not numerically trivial.}
\end{enumerate}
\end{lem}

\begin{proof} 
Consider the  decomposition 
$$W=V_1 \cup \cdots \cup V_q \cup V_{q+1} \cup \cdots \cup V_r,$$
  where 
$V_1,\dots,V_r$ are distinct  irreducible components of $W$, and assume that
$L|_{V_i}\not\equiv 0$ for $1\leq i\leq q$ and 
$L|_{V_j}\equiv 0$ for $q+1\leq j\leq r$. 
We may assume $q\ge 1$. 
Let $W_1:=V_2\cup \dots\cup V_r$. We claim that there exists a positive integer $\lambda_1$ and, for any $\lambda\ge \lambda_1$, there exists an ample $\mathbb R$-line bundle $A_{\lambda}^{(1)}$ such that 
\begin{itemize}
\item{$\frac{1}{4}A-A_{\lambda}^{(1)}$ is ample,}
\item{$\lambda L+A_{\lambda}^{(1)}$ is a $\mathbb Q$-line bundle, and}
\item{$H^0(V_1,l(\lambda L+A_{\lambda}^{(1)})\otimes (m_x^{al}   \cap I_{V_1\cap W_1}   ))\neq 0$ 
for any sufficiently divisible  integer $l>0$. }
\end{itemize}
We may assume $x\in V_1$ otherwise the result is obvious. 
Hence, Lemma \ref{l-cut1} implies the claim. 

Fix $\lambda\geq \lambda_1$. 
Then, for any sufficiently divisible $l>0$, we can find a  non-zero section  
$$0\neq s_1\in H^0(V_1,l(\lambda L+A_{\lambda}^{(1)})\otimes (m_x^{al}   \cap I_{V_1\cap W_1}   )).$$
After possibly replacing $l$ by its multiple, we can find  $t'_1\in H^0(W,l(\lambda L+A_{\lambda}^{(1)}))$ such that $t'_{1}|_{V_1}=s_1$ and $t'_{1}|_{W_1}=0$. 
In particular, 
$$t'_1\in H^0(W,l(\lambda L+A_{\lambda}^{(1)}) \otimes (m_x^{al}+I_{V_1}) ).$$
Let $B$ be an ample $\mathbb Q$-line bundle such that $\frac 1 4 A-B$ is ample. 
We may assume that $l$ is sufficiently large so that there exists 
$$t''_1\in H^0(W, lB\otimes I_{W_1})$$
such that $t''_{1}|_{V_1}\neq 0$. Let $t_1=t'_1t''_1$. Then $t_1\in H^0(W,l(\lambda L+A_{\lambda}^{(1)}+B)\otimes m^{al}_x)$ is such that 
$$t_{1}|_{V_1}\neq 0\quad\text{and}\quad t_{1}|_{V_j}=0\quad\text{for }j\neq 1.$$
We define $D_{\lambda}^{(1)}:=A_{\lambda}^{(1)}+B.$  Note that the $\mathbb R$-line bundles $D_{\lambda}^{(1)}$ and $\frac 1 2 A-D_{\lambda}^{(1)}$ are ample for any $\lambda\ge \lambda_1$. 

Similarly, we can find positive integers $\lambda_2,\dots,\lambda_q$ and sequences of 
ample $\mathbb R$-line bundles $\{D_{\lambda}^{(2)}\}_{\lambda\geq \lambda_2}, \cdots, \{D_{\lambda}^{(q)}\}_{\lambda\geq \lambda_q}$ such that if $\lambda_0:=\max\{\lambda_i\}$ 
then 
\begin{itemize}
\item{$\frac 1 2 A-D_{\lambda_0}^{(i)}$ is ample,}
\item{$\lambda_0L+D_{\lambda_0}^{(i)}$ is a $\mathbb Q$-line bundle, and}
\item{for any sufficiently divisible integer $l>0$, 
there exists $$t_i\in H^0(W, l(\lambda_0 L+D_{\lambda_0}^{(i)})\otimes m^{al}_x)$$ 
such that 
$$t_{i}|_{V_i}\neq 0\quad\text{and}\quad t_{i}|_{V_j}=0\quad\text{for }j\neq i.$$}
\end{itemize}
We define an ample $\mathbb R$-line bundle $A'$ such that 
\begin{itemize}
\item{$A-A'$ is ample,}
\item{$A'-\frac 1 2 A$ is ample, and}
\item{$\lambda_0 L+A'$ is a $\mathbb Q$-line bundle.}
\end{itemize}
Then, 
$$A'-D_{\lambda_0}^{(i)}=(\lambda_0 L+A')-(\lambda_0L+D_{\lambda_0}^{(i)})$$
is a $\mathbb Q$-line bundle. 
Moreover, 
$$A'-D_{\lambda_0}^{(i)}=(A'-\frac 1 2 A)+(\frac 1 2 A-D_{\lambda_0}^{(i)})$$ is ample. 
Thus, for sufficiently divisible integer $l>0$ and for any $i=1,\dots,q$, 
there exists $\overline{t}_i\in H^0(W, l(A'-D_{\lambda_0}^{(i)}))$ such that $\overline{t}_i|_{V_i}\neq 0$. 
Let 
$$u_i:=t_i\overline{t}_i\in H^0(W, l(\lambda_0 L+A')\otimes m^{al}_x).$$ 
Then, $u_i$ satisfies 
$$u_{i}|_{V_i}\neq 0\quad\text{and}\quad u_{i}|_{V_j}=0\quad\text{for }j\neq i.$$
We define $s:=u_1+\dots+u_q\in H^0(W,l(\lambda_0 L+A')\otimes m^{al}_x).$ 
Then 
$s|_V\neq 0$ for every irreducible component $V$ of $W$ such that 
$L|_V$ is not numerically trivial.
\end{proof}

\medskip 

We can now proceed with the proof of Proposition \ref{q-cut}. 

\begin{proof}[Proof of Proposition \ref{q-cut}] 
By Lemma \ref{l-cut2}, there exists a positive integer $\lambda_0$ and an ample $\mathbb R$-ample line bundle $A'$ on $X$ such that $A-A'$ is ample, $\lambda_0L+A'$ is a $\mathbb Q$-line bundle and for any sufficiently  divisible $l>0$, there exists 
$$s'\in H^0(W^{\red},l(\lambda_0L+A')|_{W^{\red}}\otimes m_x^{al})$$
such that $s'|_V\neq 0$ for every irreducible component $V$ of $W^{\red}$ such that $L|_{V}\not\equiv 0$.

 By Serre's vanishing theorem, if $l$ is sufficiently large, then 
 $$H^1(X,l(\lambda_0L+A')\otimes I_{W^{\red}})=0,$$ thus there exists a section 
$$s\in H^0(X,l(\lambda_0L+A')\otimes (m_x^{al}+I_{W^{\red}}))$$
such that $s|_{W^{\red}}=s'$.

Let $e$ be a positive integer such that 
$$(I_{W^{\red}})^{[p^e]}\subseteq I_W. $$
Then, we have:
$$\begin{aligned}
t:=(F^e)^* s&\in H^0(X, lp^e(\lambda_0L+A')\otimes (m_x^{al}+I_{W^{\red}})^{[p^e]})\\
&=H^0(X, lp^e(\lambda_0L+A')\otimes ((m_x^{al})^{[p^e]}+I_{W^{\red}}^{[p^e]}))\\
&\subseteq H^0(X,lp^e(\lambda_0L+A')\otimes (m_x^{alp^e}+I_{W}))\\
\end{aligned}
$$
and for every irreducible component $V$ of $W^{\red}$ such that $L|_V$ is not numerically trivial, we have
$$t|_V=((F^e)^* s)|_V= ((F^e)^*(s|_V)) \neq 0.$$
Thus, the claim follows. 
\end{proof}

As corollaries of Proposition \ref{q-cut}, 
we obtain the following two assertions. 

\begin{prop}\label{p-cut}
Fix $a \in \mathbb{ N}$. Let $X$ be a projective variety. Let $A$ be an ample $\mathbb R$-line bundle on $X$ and $L$ a strictly nef $\mathbb R$-line bundle on $X$. Let $x \in X$ be a closed point and let $ W$ be a proper closed subscheme of $X$. If $\dim_x W \ge 1$, then there exists a positive integer  $\lambda_0$ and an ample $\mathbb R$-line bundle $A'$ on $X$ such that 
\begin{enumerate}
\item{$A-A'$ is ample, }
\item{$\lambda_0L+A'$ is a $\mathbb Q$-line bundle, and }
\item{for any sufficiently divisible $l>0$, there exists
$$ t \in H^0(X, l(\lambda_0L+A') \otimes (m^{al}_x +  I_W))$$
such that $t|_V\neq 0$ for every irreducible component $V$ of $W^{\red}$.}
\end{enumerate}
\end{prop}

\begin{proof}We write $W^{\red}=W'\cup W''$, where $W'$ consists of positive dimensional components and $W''$ are the isolated points of $W^{\red}$.  By assumption, $x\in W'$. Therefore, it suffices to verify the statements (1)-(3) for $W'$.
Since $L$ is strictly nef, 
the claim follows from Lemma \ref{s-nef} and Proposition \ref{q-cut}. 
\end{proof}

\begin{prop}\label{p-cut2}
Assume that $X$ is a projective variety defined over an uncountable algebraically closed field of  characteristic $p>0$. 
Fix $a \in \mathbb{ N}$. Let $A$ be an ample $\mathbb R$-line bundle on $X$ and $L$ a nef $\mathbb R$-line bundle of maximal nef dimension on $X$. 
Let $x, y \in X$ be very general points and 
let $W$ be a proper closed subscheme of $X$ such that $\dim_x W \ge 1$ and $\dim_y W \ge 1$. 

Then, there exists a positive integer $\lambda_0$ and an ample $\mathbb R$-line bundle $A'$ on $X$ such that 
\begin{enumerate}
\item{$A-A'$ is ample, }
\item{$\lambda_0L+A'$ is a $\mathbb Q$-line bundle, and }
\item{for any sufficiently divisible $l>0$, one can find 
$$t \in H^0(X, l(\lambda_0L+A') \otimes (m^{al}_x m^{al}_y+  I_W))$$
such that $t|_V\neq 0$ for every irreducible component $V$ of $W^{\red}$ such that 
$x\in V$ or $y\in V$.}
\end{enumerate}
\end{prop}

\begin{proof}
Lemma \ref{l-nef2} implies that $L|_V$ is not numerically trivial for every irreducible component $V$ of $W^{\red}$ such that $x\in V$ or $y\in V$. Thus, by Proposition \ref{q-cut}, there exist a positive integer $\lambda_0$ and an ample $\mathbb R$-line bundle $A'$ on $X$ such that $A-A'$ is ample, $\lambda_0L+A'$ is a $\mathbb Q$-line bundle and for any sufficiently divisible $l>0$, 
we can find 
$$t_1\in H^0(X, \frac{l}{2}(\lambda_0L+A')\otimes (m^{al}_x +  I_W))$$
and
$$t_2\in H^0(X, \frac{l}{2}(\lambda_0L+A')\otimes (m^{al}_y +  I_W))$$
such that $t_i|_V\neq 0$ for every irreducible component $V$ of $W^{\red}$ such that 
$x\in V$ or $y\in V$. 
Then, $t:=t_1t_2$ is a required section. 
\end{proof}

\medskip

\subsection{Induction}\label{s-induction}

In this subsection, we describe an inductive method to construct a zero-dimensional subscheme from which we can lift sections. The subscheme is  obtained by taking the intersection of a sequence of suitable divisors. 

\begin{notation} 
Through this section, we assume that $X$ is a normal variety defined over an algebraically closed field of characteristic $p>0$. 
Let $B$ be an effective $\mathbb Q$-divisor such that
$K_X+B$ is a  $\mathbb Q$-Cartier $\mathbb Q$-divisor whose index is not divisible by $p$.
Assume that $(X, B)$ is sharply $F$-pure at a closed point $x \in X$.

Let $M\subseteq\mathbb{N}$ be the subset of positive integers $e$ such  $(p^e-1)(K_X+B)$ is Cartier. For any $i=1,\dots,r$, let $t_i:M\to \mathbb{Z}_{\ge 0} $ be a function and let $D_i$ be an effective divisor on $X$. Let $M'\subseteq M$ be an infinite subset. By abuse of notation, we say that the pair $(X,B+\sum t_iD_i)$ is {\it $M'$-sharply $F$-pure} at a point $x\in X$ if the trace map
$$\Tr^e : F^e_*(\MO_X(-(p^e - 1)(K_X+B) -\sum^r_{i=1} t_i(e)D_i)) \to  \MO_X$$
is surjective locally around $x$ for every $e\in M'$. 

Assume that the pair $(X,B+\sum^r_{i=1}t_iD_i)$ is $M'$-sharply $F$-pure at $x\in X$. Let $D_{r+1}$ be an effective divisor on $X$ such that $x\in \Supp D_{r+1}$. Then for any $e\in M'$, we  denote by 
$$\nu^{m_x}_{p^e}(X,B+\sum^r_{i=1}t_i(e)D_i;D_{r+1})$$ 
 the $F$-{\em threshold function}  of $(X,B+\sum^r_{i=1}t_i(e)D_i)$  at $x$ with respect to $D_{r+1}$, which is the maximum integer $t\ge 0$ such that the trace map  
$$\Tr^e : F^e_*(\MO_X(-(p^e - 1)(K_X+B) -\sum^r_{i=1} t_i(e)D_i-tD_{r+1})) \to  \MO_X$$
is surjective locally around $x$ (see \cite{MTW05}).

\end{notation}

\begin{say}   
Let $(X,B)$ be an $n$-dimensional  sharply $F$-pure projective pair such that $B$ is an effective ${\mathbb Q}$-divisor. Assume that the index of $K_X+B$ is not divisible by $p$.  Let $M\subseteq \mathbb{N}$  be the subset   of  positive integers $e$ such that $(p^e-1)(K_X+B)$ is Cartier. Let $A$ be an ample $\mathbb R$-Cartier $\mathbb R$-divisor on $X$ and let $L$ be a strictly nef $\mathbb R$-Cartier $\mathbb R$-divisor on $X$. Fix $a\in \mathbb{N}$. Let 
$$n_0:=\max\{\dim_k(m_x/m^2_x)\mid x \mbox{ is a closed point of $X$}\}$$
to be the maximal embedding dimension of $x\in X$.  Pick a closed point $x\in X$. Fix an integer $0\le r<n$. 
 
We assume that we have quintuples $(l_i,\lambda_i, t_i, D_i, A_i)$ for  $0\le i\le r$ where $l_i$ and $\lambda_i$ are positive integers, $t_i:M\to \mathbb{Z}_{\ge 0}$ is a function, $D_i$ is an effective Cartier divisor on $X$ and  $A_i$ is an ample $\mathbb R$-Cartier  $\mathbb R$-divisor. We assume that if $i=0$ then
$$(l_0,\lambda_0, t_0, D_0, A_0):=(0,0,0,0, A),$$ and for $i=1,\dots,r$, the quintuple $(l_i,\lambda_i, t_i, D_i, A_i)$
  satisfies the following properties:
 \begin{enumerate}
\item[$(1)_{r} $] $A-A_i$ is ample for every $1 \le i \le r$,
\item[$(2)_{r} $]  $\lambda_iL+A_i$ is a $\mathbb Q$-Cartier $\mathbb Q$-divisor, 
$l_i(\lambda_iL+A_i)$ is Cartier and 
$l_i(\lambda_iL+A_i) \sim D_i$ for every $1 \le i \le r$,
 \item[$(3)_{r}$]    $(X, B+ \sum^r_{i=1} t_iD_i)$ is $M$-sharply $F$-pure at $x$,
\item[$(4)_{r}$]  $x\in W_r$ where $W_r:=\bigcap^r_{i=1}D_i$,
 \item[$(5)_{r}$]   $\dim_xW_r = n -r$,
\item[$(6)_{r}$]  $0\le t_i(e)< \lceil\frac{n_0p^e}{al_i} \rceil $ for every $1\le i\le r$ and $e\in M$, and
\item[$(7)_{r}$]  assuming that $\Tr^e : F^e_*\MO_X(-(p^e - 1)(K_X+B) - \sum_{i=1}^r t_i(e)D_i) \to  \MO_X$ is the trace map, 
we have
$$\Tr^e (F^e_*(\MO_X(-(p^e - 1)(K_X+B) - \sum_{i=1}^r  t_i(e)D_i)\cdot I_{W_r})) \subseteq m_x,$$
for any $e\in M$. 
\end{enumerate}

We now want to construct a quintuple $(l_{r+1},\lambda_{r+1},t_{r+1},D_{r+1},A_{r+1})$, so that for 
$i=1,\dots,r+1$, the quintuple 
 $(l_i,\lambda_i, t_i, D_i, A_i)$ satisfies the above properties (1)$_{r+1}$-(7)$_{r+1}$. 

To this end, note that  Proposition \ref{p-cut} implies that there exist 
a positive integer $\lambda_{r+1}$, an ample $\mathbb R$-Cartier $\mathbb R$-divisor $A_{r+1}$ and a sufficiently divisible  integer $l_{r+1}>0$
such that 
\begin{itemize}
\item{$A-A_{r+1}$ is ample,}
\item{$\lambda_{r+1}L+A_{r+1}$ is a $\mathbb Q$-Cartier $\mathbb Q$-divisor,  $l_{r+1}(\lambda_{r+1}L+A_{r+1})$ is Cartier, and}
\item{there exists
$$s \in H^0(X, l_{r+1}(\lambda_{r+1}L+A_{r+1}) \otimes (m^{al_{r+1}}_x +  I_{W_r}))$$
such that $s|_V\neq 0$ for every irreducible component $V$ of $W_r^{\red}$. }
\end{itemize}
Let $D_{r+1}$ be the effective Cartier divisor on $X$ corresponding to $s$ and for any $e\in M$ let
$$t_{r+1}(e):=\nu^{m_x}_{p^e}(X,B+\sum^r_{i=1}t_i(e)D_i;D_{r+1}).$$ 

We now check that the properties  (1)$_{r+1}$-(7)$_{r+1}$ hold.  First note that  (1)$_{r+1}$-(5)$_{r+1}$ hold simply by  the assumptions above. 
We now check 
(6)$_{r+1}$. It is sufficient to show that for any $e\in M$ the trace map 
$$\Tr^e : F^e_*( \MO_X(-(p^e - 1)(K_X+B) -\sum^r_{i=1} t_i(e)D_i -\lceil  \frac{n_0p^e}{al_{r+1}} \rceil D_{r+1}) )\to \MO_X$$
is not surjective locally around $x$. We take an affine open subset $x\in \spec R \subseteq X$ such that $m_x|_{\spec R}$ is generated by at most $n_0$ elements and that
$$\MO_X(-D_{r+1})|_{\spec R} = fR.$$ By the definition of $D_{r+1}$, we can write
$f=\mu+\nu$
where $\mu \in m_x^{al_{r+1}}$ and $\nu \in  I_{W_r} |_{\spec R}$. Thus,  
$$f ^{\lceil  \frac{n_0p^e}{al_{r+1}} \rceil } = \mu^{\lceil  \frac{n_0p^e}{al_{r+1}} \rceil } + \nu',$$
with $\nu'\in I_{W_r} |_{\spec R}$. Then, we have
$$\begin{aligned}
\MO_X(-\lceil  \frac{n_0p^e}{al_{r+1}} \rceil D_{r+1})|_{\spec R}&=f ^{\lceil  \frac{n_0p^e}{al_{r+1}} \rceil }R\\
&=(\mu^{\lceil  \frac{n_0p^e}{al_{r+1}} \rceil } + \nu')R\\
&\subseteq m_x^{al_{r+1}   \lceil  \frac{n_0p^e}{al_{r+1}} \rceil   }+I_{W_r}|_{\spec R}\\
&\subseteq m_x^{n_0p^e}+I_{W_r}|_{\spec R}\\
&\subseteq m_x^{[p^e]}+I_{W_r}|_{\spec R} .
\end{aligned}
$$
The last inclusion was obtained as a consequence of the fact that by assumption $m_x|_{\spec R}$ is generated by at most $n_0$ elements.
We claim that  $\Tr^e(m_x^{[p^e]})\subseteq m_x$. If $f \in m_x^{[p^e]}$, then $V(f)\ge p^eD$ for some effective divisor $D$ in a neighborhood of $x$ with $x\in {\rm Supp}(D)$, which implies that $V(\Tr^e(f))\ge D$. Thus, $\Tr^e(f)\in m_x$, as claimed. It follows that 
$$\begin{aligned}
& \Tr^e \left( F^e_*\Big( \MO_X(-(p^e - 1)(K_X+B) -\sum^r_{i=1} t_i(e)D_i -\lceil  \frac{n_0p^e}{al_{r+1}} \rceil D_{r+1}\Big) \right)\\
\subseteq& \Tr^e \left( F^e_*\Big( \MO_X(-(p^e - 1)(K_X+B) -\sum^r_{i=1} t_i(e)D_i )(m_x^{[p^e]}+I_{W_r}\Big) \right)\\
\subseteq&\Tr^e(m_x^{[p^e]}) +\Tr^e \left( F^e_*\Big( \MO_X(-(p^e - 1)(K_X+B) -\sum^r_{i=1} t_i(e)D_i )I_{W_r}\Big) \right)\\
\subseteq &m_x+m_x=m_x,
\end{aligned}
$$
where we just proved the first inclusion and the third inclusion follows from (7)$_r$. 

We now prove (7)$_{r+1}$. By construction, we have $I_{W_{r+1}}=I_{W_r}+I_{D_{r+1}}$. By definition of $\nu^{m_x}_{p^e}(X,B+\sum^r_{i=1}t_i(e)D_e;D_{r+1})$, we have that
$$\Tr^e (F^e_*\MO_X(-(p^e - 1)(K_X+B) -\sum^{r+1}_{i=1} t_i(e)D_i - D_{r+1} )) \subseteq m_x$$
   for any $e\in M$. Thus, the claim follows.

\end{say}

\medskip 
To summarise, we have obtained the following theorem. 

\begin{thm}\label{t-cutcenter}
 Let $(X,B)$ be an $n$-dimensional  projective sharply $F$-pure pair such that the Cartier index of $K_X+B$ is not divisible by $p$. Fix $a\in \mathbb N$. 
Assume that $L$ is a strictly nef $\mathbb R$-Cartier $\mathbb R$-divisor, $A$ is an ample $\mathbb R$-Cartier $\mathbb R$-divisor and $M\subseteq N$ is the subset of positive integers $e$ such that $(1-p^e)(K_X+B)$ is Cartier.
 Let
$$n_0:=\max\{\dim_k(m_x/m^2_x)\mid x \mbox{ is a closed point of $X$}\}.$$
Fix a closed point $x\in X$.

Then, for any $1\le i\le n$, there are positive integers $l_i$ and $\lambda_i$,  an effective Cartier divisor $D_i$, 
an ample $\mathbb R$-Cartier $\mathbb R$-divisor $A_i$ 
and a function  $t_i:M\to \mathbb{Z}_{\ge 0}$ such that if we write $W=\bigcap^n_{i=1} D_i$, $D^{(e)}=\sum^n_{i=1} t_i(e)D_i$
and 
$$\mathcal{L}^{(e)}=\MO_X((1-p^e)(K_X+B)-D^{(e)})$$ 
then
\begin{enumerate}
\item $A-A_i$ is ample for every $1\le i \le n$,
\item $\lambda_iL+A_i$ is a $\mathbb Q$-Cartier $\mathbb Q$-divisor, 
$l_i(\lambda_iL+A_i)$ is Cartier and 
$l_i(\lambda_iL+A_i) \sim D_i$ for every $1 \le i \le n$,
\item $(X, B+\sum_{i=1}^nt_iD_i)$ is $M$-sharply $F$-pure at $x$,
\item $x\in W$,
\item $\dim_xW=0$,
\item $0\le t_i(e)<\frac{n_0p^e}{l_ia}$, for every $1\le i\le n$, and 
\item for any  
$e\in M$, we have ${\rm Tr}^e(F^e_*(\mathcal{L}^{(e)}\cdot I_W))\subseteq m_x$
and  there is an exact sequence 
\begin{diagram}
0&\rTo&F^e_*(\mathcal{L}^{(e)}\otimes I_W)&\rTo &F^e_*(\mathcal{L}^{(e)})&\rTo& F^e_*(\mathcal{L}^{(e)}\otimes \mathcal{O}_W)&\rTo& 0.
\end{diagram}  
\end{enumerate}

\end{thm}



\section{Proof of Theorem \ref{t-main}}\label{s-extend}
In this section we prove Theorem \ref{t-main}. We begin with the following:

\begin{lem}\label{perturb}
Let $(X, B)$ be a projective strongly $F$-regular pair such that $B$ is an effective $\mathbb R$-divisor. 
Let $A$ be an ample $\mathbb R$-Cartier $\mathbb R$-divisor. 
Let $L:=K_X+A+B$. 
Then, there exists an effective $\mathbb Q$-divisor $B'$ and 
an ample $\mathbb R$-Cartier $\mathbb R$-divisor $A'$ such that 
\begin{enumerate}
\item $(X, B')$ is strongly $F$-regular,
\item $K_X+B'$ is a $\mathbb Q$-Cartier divisor 
whose index is not divisible by $p$, and
\item{$L=K_X+A'+B'$.}
\end{enumerate}
\end{lem}

\begin{proof}
Let $E\ge 0$ be a divisor such that $E-K_X$ is Cartier and let $H=\Supp (E+B)$. Let $\varepsilon>0$ be a sufficiently small number such that for any effective $\mathbb Q$-Cartier $\mathbb Q$-divisor $D\le\varepsilon H$ we have that $A-D$ is ample and $(X,B+D)$ is strongly $F$-regular (cf. (1) of Remark \ref{F-sing-rem}).  

By \cite[Lemma 2.13]{HX13}, there exists an effective $\mathbb Q$-Cartier $\mathbb Q$-divisor $D$ such that the $\mathbb Q$-Cartier index of $K_X+B+D$ is not divisible by $p$. 
Thus, if we define $B'=B+D$ and $A'=A-D$, then the claim follows. 
\end{proof}

%

\medskip

\begin{say}

We can now proceed with the proof of our main theorem. 

\begin{proof}[Proof of Theorem \ref{t-main}] 
Let $L:=K_X+A+B$. 
Thus, $L$ is a strictly nef $\mathbb R$-Cartier $\mathbb R$-divisor.  
By Lemma~\ref{perturb}, we may assume that 
$B$ is a $\mathbb Q$-divisor and that the Cartier index of $K_X+B$ is not 
divisible by $p$. Let $n=\dim X$. 
We fix a positive integer $a$ such that $a\ge nn_0$,
where
$$n_0:=\sup\{\dim_k(m_x/m^2_x)\mid x \mbox{ is a closed point of $X$}\}.$$ 

Fix a closed point $x\in X$. 
We claim that 
$$L\sim_{\mathbb R} \sum_{j=1}^q c_jL_j$$
where $c_j\in\mathbb R_{\geq 0}$ and each $L_j$ is an effective Cartier divisor such that $x\not\in \Supp L_i$. 
By Lemma \ref{l-nef3}, the claim implies the Theorem. 

We apply Theorem \ref{t-cutcenter}. 
Then, for any $1\le i\le n$, 
we obtain positive integers $l_i$ and $\lambda_i$,  an effective Cartier divisor $D_i$, 
an ample $\mathbb R$-Cartier $\mathbb R$-divisor $A_i$ 
and a function  $t_i:M\to \mathbb{Z}_{\ge 0}$ such that if we write $W=\bigcap^n_{i=1} D_i$, $D^{(e)}=\sum^n_{i=1} t_i(e)D_i$,
$$L^{(e)}=(1-p^e)(K_X+B)-D^{(e)}$$ 
and $\mathcal{L}^{(e)}=\MO_X(L^{(e)})$,
then
\begin{enumerate}
\item $\frac{1}{2}A-A_i$ is ample for every $1\le i \le n$,
\item $\lambda_iL+A_i$ is a $\mathbb Q$-Cartier $\mathbb Q$-divisor, 
$l_i(\lambda_iL+A_i)$ is Cartier and 
$l_i(\lambda_iL+A_i) \sim D_i$ for every $1 \le i \le n$,
\item $(X, B+\sum_{i=1}^nt_iD_i)$ is $M$-sharply $F$-pure at $x$,
\item $x\in W$,
\item $\dim_xW=0$,
\item $0\le t_i(e)<\frac{n_0p^e}{l_i a}$, for every $1\le i\le n$ and $e\in M$, and
\item ${\rm Tr}^e(F^e_*(\mathcal{L}^{(e)}\cdot I_W))\subseteq m_x$, for every $e\in M$.
\end{enumerate}
In particular, we have that $D^{(e)}=\sum^n_{i=1} t_i(e)D_i\sim \sum^n_{i=1} t_i(e)l_i(\lambda_iL+A_i).$

We can write $L=\sum_{i=1}^r \alpha_iE_i$ where $\alpha_i\in\mathbb R$ and $E_i$ are Cartier divisors, for $i=1,\dots,r$. 
Let $V\subseteq {\rm Div}_{\mathbb R}(X)$ be the vector space spanned by $E_1,\dots,E_r$. 
We denote by $\|\cdot \|$ the sup norm with respect to this basis. 
Let $\varepsilon>0$ be a sufficiently small 
rational number such that $\frac 1 2 A-\Gamma$ is ample for any $\Gamma\in V$ such that $\|\Gamma\|<\varepsilon$.
By Diophantine approximation (e.g. see \cite[Lemma 3.7.7]{bchm}), there exist $\mathbb Q$-divisors $C'_j\in V$ with $j=1,\dots,q$ and 
positive integers $m_j>1+\sum^n_{i=1} \frac{n_0\lambda_i}{a}$ such that $L=\sum_{j=1}^{q}r_jC_j'$ for some 
real numbers $0\leq r_j\leq 1$ with $\sum_{j=1}^{q}r_j=1$, the divisor $C_j:=m_jC'_j$ is Cartier and 
$$\|L - C'_j\|\le \frac{\varepsilon}{m_j}\qquad \text{for any }j=1,\dots,q.$$
Let $\Gamma_j:=m_jL-C_j$ for $j=1,\dots,q$. 
Then, we obtain  
\begin{itemize}
\item{$L=\sum_{j=1}^{q} \frac{r_j}{m_j}(m_jL-\Gamma_j)$,}
\item{$\frac{1}{2}A-\Gamma_j$ is ample,}
\item{$C_j=m_jL-\Gamma_j$ is Cartier, and }
\item{$m_j>1+\sum^n_{i=1} \frac{n_0\lambda_i}{a}$.}
\end{itemize}

We want to show that $x$ is not a base point of the linear system $|C_j|$ for  $j=1,\dots,q$.  
Fix $1\leq j\leq q$. 
For all $e\in M$, we have the following diagram: 
$$\begin{CD}
0@>>>F^e_*(\mathcal{L}^{(e)}\otimes I_W)@>>>F^e_*(\mathcal{L}^{(e)})@>>> 
F^e_*(\mathcal{L}^{(e)}\otimes \mathcal{O}_W)@>>> 0\\
@. @VVV @VV{\Tr^e}V @VV{\varphi^e}V\\
0@>>>m_x@>>> \mathcal{O}_X@>>>k_x @>>> 0
\end{CD}$$
where the first vertical arrows is the inclusion given by $(7)$ above and the third vertical arrow is the natural map obtained by diagram chasing. 

Since $(X,B+\sum_{i=1}^n t_i D_i)$ is $M$-sharply $F$-pure at $x$, it follows that $$\varphi^e\colon F^e_*(\mathcal{L}^{(e)}\otimes \mathcal{O}_W) \to k_x$$ is surjective in a neighbourhood of $x$ for all $e\in M$. 
Tensoring by $\mathcal O_X(C_j)$ and taking  cohomology, we obtain 
$$\begin{CD}
H^0(X, F^e_*(\mathcal{L}^{(e)})\otimes \mathcal O_X(C_j))@>>> 
H^0(W, F^e_*(\mathcal{L}^{(e)})\otimes\mathcal O_W(C_j))\\
@VVV @VV{H^0(\varphi^e)}V\\
H^0(X, \MO_X (C_j))@>\rho>>H^0(x, \MO_X(C_j)\otimes k_x).
\end{CD}$$
Since $\dim_x W=0$, the map  $H^0(\varphi^e)$ is surjective. 
Thus, to show that $\rho$ is surjective, 
it is enough to prove $H^1(X, \mathcal{L}^{(e)}\otimes \MO_X(p^eC_j)\otimes I_W)=0.$
We have
$$\begin{aligned}
{L}^{(e)}+p^eC_j~=&~(1-p^e)(K_X+B)-\sum^n_{i=1} t_i(e)D_i+p^em_jL-p^e\Gamma_j\\
\sim&~ (1-p^e)(K_X+B)-\sum^n_{i=1} t_i(e)l_i(\lambda_iL+A_i)+p^em_jL-p^e\Gamma_j\\
=&~K_X+B+p^eA-\sum^n_{i=1} t_i(e)l_iA_i-p^e\Gamma_j+(p^em_j-p^e-\sum^n_{i=1} t_i(e)l_i\lambda_i)L\\
=&~K_X+B+(p^e-\frac 1 2 \sum^n_{i=1} t_i(e)l_i-\frac 1 2 p^e)A\\
&~\quad +\sum^n_{i=1} t_i(e)l_i(\frac 1 2 A-A_i)+p^e(\frac 1 2 A-\Gamma_j)+(p^em_j-p^e-\sum^n_{i=1} t_i(e)l_i\lambda_i)L.
\end{aligned}$$
As $a\ge nn_0$, it follows 
$$p^e-\frac 1 2 \sum^n_{i=1} t_i(e)l_i-\frac{p^e}{2}\geq\frac{p^e}{2}-\frac 1 2\sum^n_{i=1}\frac{n_0p^e}{a} \ge\frac{p^e}{2}-\frac{p^e}{2}=0.$$
Since $m_j>1+\sum^n_{i=1} \frac{n_0\lambda_i}{a}$, we have 
\begin{eqnarray*}
p^em_j-p^e-\sum^n_{i=1} t_i(e)l_i\lambda_i
\geq p^em_j-p^e-\sum^n_{i=1} \frac{n_0p^e}{a}\lambda_i> 0.
\end{eqnarray*}
Since $\frac 1 2 A-A_i$ and $\frac 1 2 A-\Gamma_j$ are ample,  the Fujita vanishing theorem implies that if $e\in M$ is  sufficiently large then 
$$H^1(X, \mathcal{L}^{(e)}\otimes \MO_X(p^eC_j)\otimes I_W)=0.$$
Thus, the claim follows. 
\end{proof}	

\medskip 

\begin{proof}[Proof of Corollary \ref{c-main}]Since $K_X+\Delta$ is big, we can write $K_X+\Delta\sim_{\mathbb{R}}A+B$ where $A$ is an ample effective $\mathbb R$-Cartier $\mathbb R$-divisor  and $B$ is effective. Replacing $\Delta$ by $\Delta'=\Delta+tB$, and choosing $t$ to be a sufficiently small positive number such that $(X,\Delta')$ is strongly $F$-regular, then the assertion follows from Theorem \ref{t-main}. 
\end{proof}
\end{say}

\medskip

\begin{say}
As an immediate consequence of Theorem \ref{t-main} we obtain the rationality theorem:

\begin{proof}[Proof of Theorem \ref{t-rational}]
Since $K_X+B$ is not nef, we have that $\lambda>0$. By the definition of $\lambda$ it follows that $K_X+\lambda A+B$ is nef but not ample. Thus,  Theorem \ref{t-main} implies that $K_X+\lambda A+B$ is not strictly nef. In particular, there exists a curve $C$ such that 
$(K_X+\lambda A+B)\cdot C=0$, i.e. 
$$\lambda=  \frac{-(K_X+B)\cdot C}{A\cdot C}.$$
Thus, $\lambda$ is rational. 
\end{proof}
\end{say}

\medskip 

\section{Proof of Theorem \ref{t-main2}}\label{s-max}

 We now proceed with the proof of Theorem \ref{t-main2}. We use the same methods we developed in the last section,  but this time we cut at two points at the same time. 
\begin{say}
\begin{proof}[Proof of Theorem \ref{t-main2}] Since it suffices to prove the statement of the theorem after any base change of 
the ground field, we may assume the ground field is uncountable. Fix  two very general points $x,y\in X$ as in Lemma 
\ref{l-nef2}. In addition, we assume that $x$ and $y$ are not contained in the singular locus of $X$, nor in the support of $B$. In particular, $m_x$ and $m_y$ are generated by $n$ elements.

By using the same argument as in the proof of Lemma~\ref{perturb}, 
we may assume that $B$ is an effective ${\mathbb Q}$-divisor such that 
the Cartier index of $K_X+B$ is not divisible by $p$. Let $n=\dim X$. We fix an integer $a>2n^2$. 
Let $M_0\subseteq\mathbb{N}$ be the subset of positive integers $e$ such  $(p^e-1)(K_X+B)$ is Cartier and let  $L:=K_X+A+B$. By assumption, $(X,B)$ is $M_0$-sharply $F$-pure at $x$ and $y$. 
By Lemma~\ref{l-nef4}, 
it is sufficient to show that 
$L\sim_{\mathbb R}\sum_{j=1}^qc_j E_j$ where $c_j\in\mathbb R_{>0}$ and 
$E_j$ is an effective Cartier divisor such that 
$\Supp E_j$ contains $x$ but not $y$ for every $1\leq j\leq q$.

We start with the seven-tuple
$$(l_0,\lambda_0, t_0, D_0, A_0,M_0,W_0):=(0,0, 0, 0, A, M_0,X).$$
Fix $0\le r < n$. Let us assume we have constructed a seven-tuple 
$(l_i,\lambda_i, t_i, D_i,A_i,M_i,W_i)$ for $0\le i\le r$ where $l_i$ and $\lambda_i$ are positive integers, 
$M_i\subseteq \mathbb{N}$ is an infinite subset,
$t_i:M_i\to \mathbb{Z}_{\ge 0}$ is a function, 
$D_i$ is an effective Cartier divisor on $X$, 
$A_i$ is an ample $\mathbb R$-Cartier $\mathbb R$-divisor on $X$, and 
$W_i$ is a closed subscheme of $X$ which satisfy either the following properties:
\begin{enumerate}
\item[(1)$'_r$]  $\frac{1}{2}A-A_i$ is ample for every $1 \le i \le r$,
\item[(2)$'_r$]  $\lambda_iL+A_i$ is a $\mathbb Q$-Cartier $\mathbb Q$-divisor, 
$l_i(\lambda_iL+A_i)$ is Cartier and 
$l_i(\lambda_iL+A_i) \sim D_i$ for every $1 \le i \le r$. 
\item[(3)$'_r$]   $(X, B+ \sum^r_{i=1} t_iD_i)$ is $M_r$-sharply $F$-pure at $x$ and $y$,
\item[(4)$'_r$] $x,y\in W_r$ where $W_r=\bigcap^r_{i=1}D_i$,
\item[(5)$'_r$]  $\dim_xW_r =\dim_y W_r= n -r$,
\item[(6)$'_r$] $0\le t_i(e)< \lceil \frac{np^e}{al_i} \rceil$ for every $1\le i\le r$, and
\item[(7)$'_r$] 
assuming that $\Tr^e : F^e_*(\MO_X(-(p^e - 1)(K_X+B) - \sum_{i=1}^r t_i(e)D_i)) \to  \MO_X$ is the trace map, 
we have
$$\Tr^e : F^e_*(\MO_X(-(p^e - 1)(K_X+B) - \sum_{i=1}^r  t_i(e)D_i)\cdot I_{W_r}) \subseteq m_x\cap m_y.$$
\end{enumerate}
or the following properties: 
\begin{enumerate}
\item[(1)$''_r$]  $\frac{1}{2}A-A_i$ is ample for every $1 \le i \le r$,
\item[(2)$''_r$]  $\lambda_iL+A_i$ is a $\mathbb Q$-Cartier $\mathbb Q$-divisor, 
$l_i(\lambda_iL+A_i)$ is Cartier and 
$l_i(\lambda_iL+A_i) \sim D_i$ for every $1 \le i \le r$,
\item[(3)$''_r$]     $(X, B+ \sum^r_{i=1} t_iD_i)$ is $M_r$-sharply $F$-pure at $x$ and 
 $(X, B+ \sum^r_{i=1} t_i(e)D_i)$
 is not $M_r$-sharply $F$-pure at $y$ for every $e\in M_r$,
\item[(4)$''_r$]   $x\in W_r$ where $W_r=\bigcap^r_{i=1}D_i$,
\item[(5)$''_r$]    $\dim_xW_r = n -r$,
\item[(6)$''_r$]   $0\le t_i(e)<\lceil \frac{np^e}{al_i} \rceil$ for every $1\le i\le r$, and 
\item[(7)$''_r$]   
assuming that $\Tr^e : F^e_*(\MO_X(-(p^e - 1)(K_X+B) -\sum_{i=1}^r  t_i(e)D_i)) \to  \MO_X$ is the trace map, 
we have
$$\Tr^e : F^e_*(\MO_X(-(p^e - 1)(K_X+B) - \sum_{i=1}^r t_i(e)D_i)\cdot I_{W_r}) \subseteq m_x\cap m_y.$$
\end{enumerate}
We claim that after possibly switching $x$ and $y$, we can find a seven-tuple 
$$(l_{r+1},\lambda_{r+1}, t_{r+1}, D_{r+1},A_{r+1},M_{r+1},W_{r+1})$$ such that either $(1)'_{r+1}-(7)'_{r+1}$ or $(1)''_{r+1}-(7)''_{r+1}$ hold.

Let us prove the claim.  Assume first that the properties  $(1)'_r-(7)'_r$ hold. 
Then, by Proposition \ref{p-cut2}, there exist 
positive integers $l_{r+1}$ and $\lambda_{r+1}$ and an ample  $\mathbb R$-Cartier $\mathbb R$-divisor $A_{r+1}$ on $X$ 
such that 
\begin{itemize}
\item {$\frac{1}{2}A-A_{r+1}$ is ample, }
\item {$\lambda_{r+1}L+A_{r+1}$ is a $\mathbb Q$-Cartier $\mathbb Q$-divisor,}
\item {$l_{r+1}(\lambda_{r+1}L+A_{r+1})$ is Cartier, and}
\item {there exists a section 
$$s\in H^0(X, l_{r+1}(\lambda_{r+1}L+A_{r+1})\otimes (m^{al_{r+1}}_{x}m_y^{al_{r+1}}+I_{W_r})) $$
such that $s|_V\neq 0$ for every irreducible component $V$ of $W^{\red}_r$ such that $x\in V$ or $y\in V.$}
\end{itemize}

Let $D_{r+1}$ be the effective divisor on $X$ corresponding to $s$. 
We define $W_{r+1}=W_r\cap D_{r+1}$. 
Let $$t^x_{r+1}(e):=\nu^{m_x}_{p^e}(X,B+\sum^r_{i=1}t_i(e)D_i;D_{r+1})$$ and $$t^y_{r+1}(e):=\nu^{m_y}_{p^e}(X,B+\sum^r_{i=1}t_i(e)D_i;D_{r+1})).$$ 

We  consider the sets 
\begin{itemize}
\item $M^>_{r} :=\{e\in M_r\,|\,t^x_{r+1}(e)>  ~ t^y_{r+1}(e)\}$,
\item $M^<_{r}  :=\{e\in M_r\,|\,t^x_{r+1}(e) < ~ t^y_{r+1}(e)\}$, and
\item $M^=_{r} :=\{e\in M_r\,|\,t^x_{r+1}(e)=~ t^y_{r+1}(e)\}.$
\end{itemize}

If $M_r^{=}$ is an infinite set, then we choose $M_{r+1}=M_r^=$ and $t_{r+1}(e)=t^x_{r+1}(e)$. As in the proof of Theorem \ref{t-main}, it follows that the seven-tuples $(l_{i},\lambda_{i}, t_{i}, D_{i},A_{i},M_{i},W_{i})$, with $i=1,\dots,r+1$, satisfy $(1)'_{r+1}-(7)'_{r+1}$.

If $M_r^{=}$ is not an infinite set, then one of the sets $M_r^{>}$ and $M_r^{<}$, say $M_r^{>}$, is infinite. We choose $M_{r+1}=M_r^{>}$, and $t_{r+1}(e)=t_{r+1}^x(e)$ and, as above, we easily see  that  the seven-tuples $(l_{i},\lambda_{i}, t_{i}, D_{i},A_{i},M_{i},W_{i})$, with $i=1,\dots,r+1$, satisfy $(1)''_{r+1}-(7)''_{r+1}$.

\medskip 

Now let us assume that $(1)''_r-(7)''_r$  hold. Then we ignore $y$ and just do the same construction as in the proof of Theorem \ref{t-main} for $x$, where we choose $W_{r+1}$ with the methods described above. 
To proceed with the induction, note that  Proposition \ref{p-cut} implies that there exists a positive integer $\lambda_{r+1}$, a sufficiently large and divisible  integer $l_{r+1}$ and a  section 
$$0\neq s \in H^0(X, l_{r+1}(\lambda_{r+1}L+A_{r+1}) \otimes (m^{al_{r+1}}_x +  I_{W_r}))$$
such that $s|_V\neq 0$ for every irreducible component $V$ of $W_r^{\red}$ such that 
$x\in V$. Thus, we let $D_{r+1}$ to be the corresponding effective divisor on $X$ and $t_{r+1}(e)=\nu^{m_x}_{p^e}(X,B+\sum^r_{i=1}t_i(e)D_i;D_{r+1}).$ 
In particular, $M_{r+1}=M_{r}$. Then it is easy to see that $(1)''_{r+1}-(7)''_{r+1}$ hold for the seven-tuples $(l_i,\lambda_i, t_i, D_i,M_i,W_i)$, with $i=1,\dots,r+1$. 
Thus, we have proven the claim. 

\medskip

We now apply the same argument as in Theorem \ref{t-cutcenter}. 
For any $1\le i\le n$, 
we obtain a quintuple $(l_i(x), \lambda_i(x), D_i(x), A_i(x), t_i(x))$
where 
$l_i(x)$ and $\lambda_i(x)$ are positive integers,  
$D_i(x)$ is an effective Cartier divisor, 
$A_i(x)$ is an ample $\mathbb R$-Cartier $\mathbb R$-divisor,  and 
$t_i(x):M_0\to \mathbb{Z}_{\ge 0}$ is a function such that if we write 
$W(x)=\bigcap^n_{i=1} D_i(x)$, $D^e(x)=\sum^n_{i=1} t_i(x)(e)D_i(x)$
and 
$$\mathcal{L}^e(x)=\MO_X((1-p^e)(K_X+B)-D^{e}(x))$$ 
then
\begin{enumerate}
\item[$(1)_x$] $\frac{1}{2}A-A_i(x)$ is ample for every $1\le i \le n$,
\item[$(2)_x$] $\lambda_i(x)L+A_i(x)$ is a $\mathbb Q$-Cartier $\mathbb Q$-divisor, 
$l_i(x)(\lambda_i(x)L+A_i(x))$ is Cartier and 
$l_i(x)(\lambda_i(x)L+A_i(x)) \sim D_i(x)$ for every $1 \le i \le n$,
\item[$(3)_x$] $(X, B+\sum_{i=1}^nt_i(x)D_i(x))$ is $M$-sharply $F$-pure at $x$,
\item[$(4)_x$] $x\in W(x)$,
\item[$(5)_x$] $\dim_xW(x)=0$,
\item[$(6)_x$] $0\le t_i(x)(e)<\frac{np^e}{al_i(x)}$, for every $1\le i\le n$, and
\item[$(7)_x$] ${\rm Tr}^e(F^e_*(\mathcal{L}^{e}(x)\cdot I_{W(x)}))\subseteq m_x$.
\end{enumerate}
We define a quintuple $(l_i(y), \lambda_i(y), D_i(y), A_i(y), t_i(y))$ in the same way. 

\medskip

By Diophantine approximation (see the proof of Theorem~\ref{t-main}), 
there exists  $\mathbb R$-Cartier $\mathbb R$-divisors $\Gamma_1,\dots,\Gamma_q$ 
positive integers $m_j$ and positive real numbers $r_j$ for $j=1,\dots,q$ 
such that
\begin{enumerate}
\item[($a$)]{$L=\sum_{j=1}^q \frac{r_j}{m_j}(m_jL-\Gamma_j)$,}
\item[($b$)]{$\frac{1}{2}A-\Gamma_j$ is ample,}
\item[($c$)]{$C_j:=m_jL-\Gamma_j$ is Cartier, and }
\item[($d$)]{$m_j>1+\max\{\sum^n_{i=1} \frac{n\lambda_i}{a}, \sum^n_{i=1} \frac{n\lambda_i(x)}{a}, 
\sum^n_{i=1} \frac{n\lambda_i(y)}{a}\}$.}
\end{enumerate}

It is sufficient to show that the linear system $|C_j|$ separates $x$ and $y$ for every $j=1,\dots,q$.  
Thus, we want to  prove that 
there exist sections $s,t\in H^0(X, C_j)$ such that  $s|_x\neq 0$ but $s|_y =0$ and $t|_x=0$ but $t|_y \neq 0$. 
For any $e\in M_n$, let 
$$D^{(e)}=\sum_{i=1}^n t_i(e)D_i\qquad \text{and} \qquad  \mathcal{L}^{(e)}=\MO_X((1-p^e)(K_X+B)-D^{e}).
$$ 
We first assume that $(1)'_n-(7)'_n$ is true. Then, as in the proof of Theorem \ref{t-main}, we know that for any $e\in M_n$ there is a diagram 
$$\begin{CD}
0@>>>F^e_*(\mathcal{L}^{(e)}\otimes I_W)@>>>F^e_*(\mathcal{L}^{(e)})@>>> 
F^e_*(\mathcal{L}^{(e)}\otimes \mathcal{O}_W)@>>> 0\\
@. @VVV @VV{\Tr^e}V @VV{\varphi^e}V\\
0@>>>m_x\cap m_y@>>> \mathcal{O}_X@>>>k_x\oplus k_y@>>> 0.
\end{CD}$$

By $(2)_r'$, $\varphi^e$ is surjective for any $e\in M_n$. 
Thus, as in the proof of Theorem \ref{t-main}, after tensoring by $C_j$, the diagram induces a surjection 
$$H^0(X,C_j)\to H^0(x,C_j)\oplus H^0(y,C_j).$$

On the other hand, if $(1)''_n-(7)''_n$ holds, the same diagram as above holds, but 
by $(2)_r''$, the map
$$\varphi^e\colon F^e_*(\mathcal{L}^{(e)}\otimes \mathcal{O}_{W_n})\to k_x\oplus k_y $$ factors through a map 
$$F^e_*(\mathcal{L}^{(e)}\otimes \mathcal{O}_{W_n})\to k_x\oplus 0\to k_x\oplus k_y$$ for every $e\in M_n$.
Thus, there is a section $s\in H^0(X, C_j)$ 
such that $s|_x\neq 0$ but $s|_y =0$. 
Thanks to $(1)_y-(7)_y$ and $(a)-(d)$, 
by the same argument as in Theorem \ref{t-cutcenter}, we can show that
$y$ is not a base point of  the linear system $|C_j|$. 
Then, 
we can find a section $s'$ such that $s'|_{y}\neq 0$. Hence, using a linear combination of $s$ and $s'$, we can  find a section which vanishes at $x$ but not at $y$. 

Therefore, we can apply Lemma \ref{l-nef4} to conclude that $L$ is big. 
\end{proof}
\end{say}
\medskip
\begin{say}Corollary \ref{c-rational-cover}  follows immediately from the following lemma. 
\begin{lem}\label{l-bb} 
Let $X$ be a normal projective variety, defined over an algebraically closed field $k$ of characteristic $p>0$. Assume that $A$ is an ample $\mathbb R$-divisor, $B\ge 0$ is an $\mathbb R$-divisor such that $L=K_X+A+B$ is nef but not big. Then $X$ is covered by rational curves $R$ such that 
$$L\cdot R=0 \qquad \text{and} \qquad (K_X+B)\cdot R\ge -2\dim X.$$
\end{lem}
\begin{proof} 
Let $K\supseteq k$ be an uncountable algebraically closed field. Since ${ NE}(X)={NE}(X_K)$, if there are rational curves in the class $[R]\in NE(X_K)$ covering $X_K$, then  there is a component $V$ of  ${\rm RatCurve}^n(X_K)$ parameterizing moving curves which are in $[R]$ (see \cite[Definition - Proposition II.2.11]{Kol96}). 
By \cite[II.2.15]{Kol96}, and since the construction of Hilbert schemes commutes with base change, it follows that 
$${\rm RatCurve}^n(X_K)={\rm RatCurve}^n(X)\times_kK.$$
 Thus, there exist rational curves in $R$ which cover $X$. Therefore, we may assume that the ground field is uncountable. 

Since $L$ is not big, Theorem \ref{t-main2} implies that  $L$ is not of maximal nef dimension.
Let $f\colon X\dashrightarrow Z$ be the nef reduction map associated to $L$ and whose existence is 
guaranteed by Theorem \ref{t_nrm}. Let $X'$ be the normalization of the graph $\Gamma(f)\subseteq X\times Z$. Note that the induced morphism $p_1:X'\to X$ is an isomorphism over an open set $V=f^{-1}(U)$ for some nonempty open set $U\subseteq Z$. 

Theorem \ref{t_nrm} implies that $p_1^*L$ is numerically trivial on any fibre of $p_2:X'\to Z$, i.e. $-p_1^*(K_X+B)$ is $p_2$-ample. Therefore, we can take a sufficiently ample divisor $H$ on $Z$ such that $H'=-p_1^*(K_X+B)+p_2^*H$ is ample. Furthermore, we can assume that for any curve $C$ on $X'$ which is not contained in fibres of $p_2$, we have 
$C\cdot H'> 2\dim X.$

Let $x$ be a very general point of $X'$ and let $C$ be a  curve passing through $x$ and which is contained in a fibre $F$ of $f$. We may 
assume that $C$ does not intersect the singular locus of $X'$ and it is not contained in $p_1^{-1}(\Supp B)$.  In particular, it follows that
$H'\cdot C\le -K_{X'}\cdot C$.  By Theorem \ref{t_nrm}, we have that  $L\cdot 
C=0$. 
Applying Miyaoka-Mori's bend and break (see \cite{MM86}, \cite[Theorem II.5.8]{Kol96}), it follows that there is a rational curve $R'$ passing though $x$ such that
$$H'\cdot R'\le 2\dim X \frac{H'\cdot C}{-K_{X'}\cdot C}\le2\dim X .$$
Therefore, $R'$ is contained in a fibre over $Z$. In particular, we can assume $p_1:X'\to X$ is an isomorphism on a neighbourhood of the curve $R'$ and if we denote by $R$ the image of $R'$ in $X$ then we have $R\cdot L=0$.  In addition, we have 
$$-(K_X+B)\cdot R= H'\cdot R' \le 2 \dim X$$ 
and the claim follows. 
\end{proof}
\end{say}

\section{Three dimensional MMP}
In this section, we focus on the study of three dimensional varieties defined over an algebraically closed field of positive characteristic. In Subsection \ref{s-cone}, using results from  \cite{Kol91,Keel99,HX13}, we show that a weak version of the minimal model program holds for terminal threefolds. In Subsection \ref{s-bpf}, we prove, under some restrictions on the coefficients of the boundary, that the base point free theorem holds for three dimensional log canonical pairs with intermediate Kodaira dimension. 

\subsection{A weak cone theorem and running the MMP}\label{s-cone}
The aim of this section is to prove Theorem \ref{t-cone} and Theorem \ref{t-3mmp}.

We begin with the following: 

\begin{lem}\label{l-cut} Let $X$ be a $\mathbb Q$-factorial projective variety defined over an algebraically closed field.  
Let $B$ be an effective $\mathbb R$-divisor on $X$,  let $\lambda_H$ be the  nef threshold of $K_X+B$ with respect to an ample $\mathbb R$-divisor $H$ and let 
$$\mathcal S=\{C\in N_1(X)\mid -(K_X+B)\cdot C \ge 0 \text{ and } (K_X+B+\lambda_H H )\cdot C\ge 0\text{ for all ample } H\}.$$  
Then 
$$\overline {NE}~(X)=\mathcal S\cup \overline {NE}~(X)_{K_X+B\ge 0}.$$
\end{lem}
\begin{proof}
Clearly the left hand side is contained in the right hand side. 

Assume that    
there exists $\xi$ in the interior of $\mathcal{S}$ such that $\xi\notin  \overline {NE}~(X)$.
Then there exists an hyperplane
which separates $\xi$  from $ \overline { NE}~(X)$, i.e. there exists a divisor $L$ such that $L\cdot \xi <0$ and $L\cdot C>0$ for all 
$C\in \overline { NE}~(X)\setminus \{0 \}$. In particular $L$ is ample. Since $\xi$ is contained in the interior of $\mathcal S$ it follows that
$$(K_X+B+\lambda_L L) \cdot \xi>0.$$
Thus
$$(K_X+B)\cdot \xi > -\lambda_L L\cdot \xi \ge 0,$$
which is a contradiction, since $\xi \in \mathcal S$.
\end{proof}

\begin{lem}\label{l-cone} Let $X$ be a $\mathbb Q$-factorial projective variety defined over an algebraically closed field.  Let $A$ be an ample $\mathbb{R}$-divisor on $X$ and let $B$ be an effective $\mathbb R$-divisor on $X$. For any ample $\mathbb R$-divisor $H$, let $a_H$ be the nef threshold of 
$K_X+\frac 1 2 A+B$ with respect to $H$.

Assume that there exist finitely many extremal rays of $\overline{NE}(X)$ spanned by the classes of curves $R_1$,..., $R_m$, such that for any ample $\mathbb R$-divisor $H$ on $X$, we have that 
$\overline{NE}(X)\cap (K_X+\frac{1}{2}A+B+a_HH)^{\perp}$ contains $ R_i$ for some $i$. 
Then
$$\overline{NE}(X)=\overline {NE}~(X)_{K_X+A+B\ge 0}+\sum_i \mathbb R_{\ge 0} R_i.$$
\end{lem}
\begin{proof}First we prove the left hand side is  equal to the closure of the right hand side. If not, there exists $C\in \overline{NE}(X)$ and a divisor $L$ such that 
$L\cdot \xi >0$ for any $\xi\neq 0$ in the closure of the right hand side but $L\cdot C<0$. We can assume $C$ is in the  boundary of $\overline{NE}(X)$. In particular, $L\cdot R_i>0$ for $i=1,\dots,m$.

 For any ample $\mathbb R$-divisor $H$, let $b_H$ be the nef threshold of $K_X+A+B$ with respect to $H$. 
By Lemma \ref{l-cut}, we can assume that there exists a sequence of ample divisors $H_j$ with $j\ge 1$ 
such that  $$ (K_X+A+B+b_{H_j}H_j)\cdot C<\frac{1}{j}.$$ 
Fix a sufficiently small positive number $a$ such that $\frac{1}{2}A+aL$ is ample. Fix a sufficiently large positive integer $j$, so that
$$ (K_X+A+B+b_{H_j}H_j)\cdot C<\frac{1}{j}<-aL\cdot C$$ 
and let $H'=\frac{1}{2}A+b_{H_j}H_j+aL$. Then, by assumption, the nef threshold $a_{H'}$ of $K_X+B+\frac{1}{2}A$ with respect to $H'$  is larger than 1. But 
for any $i= 1,\dots,m$, since $K_X+A+B+b_{H_j}H_j$ is nef, we have 
$$
\begin{aligned}
(K_X+\frac{1}{2}A+B+a_{H'}H')\cdot R_i
&=(K_X+A+B+b_{H_j}H_j +aL + (a_{H'}-1)H')\cdot R_i\\
&> a L\cdot R_i>0,
\end{aligned}
$$
which contradicts our assumption. 

It remains to show that $\mathcal P=\overline {NE}~(X)_{K_X+A+B\ge 0}+\sum_i \mathbb R_{\ge 0} R_i$ is closed. We use a standard argument for this.  Let $z_j\in \mathcal P$ be a sequence of points, with $j\ge 1$ such that $\lim_{j}z_j=z\in \overline{NE}(X)$. 
Then, for any $j\ge 1$, we may write  $z_j=v_j+\sum^m_{i=1} a_{ij}R_i$ for some $v_j\in \overline {NE}~(X)_{K_X+A+B\ge 0}$ and $a_{ij}\in \mathbb R_{\ge 0}$. 
Let $H$ be an ample divisor on $X$. Then intersecting with $H$, we have that if $j$ is sufficiently large, 
$H\cdot z_j\le z\cdot H+1$ and in particular it follows that the coefficients
$a_{ij}$ are bounded by a fixed constant. Thus, after passing through a subsequence, we can assume that for each $i=1,\dots,m$, the sequence $a_{ij}$ has a limit, say $a_i$. Then
$$z-\sum_{i=1}^m a_{i}R_i=\lim_i (z_i-\sum_{i=1}^m a_{ij} R_i)\in \overline{NE}(X)_{K_X+A+B\ge 0}.$$
Thus, $z\in \mathcal P$. 
\end{proof}

We now proceed with the proof of Theorem \ref{t-cone}. 
\begin{proof}[Proof of Theorem \ref{t-cone}] 
For any ample $\mathbb R$-divisor $H$, let $\lambda_H$ be the  nef threshold of $K_X+B$ with respect to $H$.

We first assume that there exists an ample $\mathbb R$-divisor $H$ such that $K_X+B+\lambda_H H$  is big.  Let $t$ be a rational number such that $0<t<\lambda_H$ and $K_X+B+t H$ is big. 
Then, by perturbing $tH$, we can find an ample $\mathbb Q$-divisor $A$ such that 
$K_X+B+A$ is big and not nef. 
Then, the result follows from \cite[Proposition 0.6]{Keel99}. 

Thus, we may assume that $K_X+B+\lambda_H H$ is not big for all ample $\mathbb R$-divisors $H$. Pick any ample $\mathbb R$-divisor $A$ such that $K_X+A+B$ is not pseudo-effective. Thus, for any ample $\mathbb R$-divisor $H$,  if $a_H$ is the nef threshold of $K_X+\frac{1}{2}A +B$ with respect to $H$, then $K_X+\frac{1}{2}A+B+ a_H H$ is not big.  
By Lemma \ref{l-bb}, there exists a rational curve $R$ such that 
$$(K_X+\frac{1}{2}A+B+a_HH)\cdot R=0, \qquad -(K_X+B)\cdot R<6, $$
which implies $A\cdot R<12$. In particular, $R$ is parametrized by finitely many components of the Chow variety ${\rm Chow_1}(X)$ and we may assume that there exists finitely many curves $R_1,\dots,R_m$ such that if $H$ is an  ample $\mathbb R$-divisor, then $(K_X+\frac{1}{2}A+B+a_HH)\cdot R_i=0$ for some $i\in \{1,\dots,m\}$. 
Thus, the result follows from Lemma \ref{l-cone}. 
\end{proof}

\medskip 

\begin{say}

We now proceed with the proof of Theorem \ref{t-3mmp}.
Case (1) of Theorem \ref{t-3mmp} is proven in  \cite{HX13}. Thus, we only need to consider the case when $K_X$ is not pseudo-effective. In this case, our result follows directly from a combination of Theorem \ref{t-cone} and Koll\'ar's contraction theorem \cite[Section 4]{Kol91}.

\begin{proof}[Proof of Theorem \ref{t-3mmp}] If $K_X$ is not nef, 
then Theorem \ref{t-cone} implies that there exists an extremal ray $R$ of $\overline{NE}(X)$ and an ample $\mathbb{Q}$-divisor $H$, such that
$K_X+H$ is nef and 
$$(K_X+H)^{\perp}\cap \overline{NE}(X)=R.$$
If $K_X+H$ is not big, then Lemma \ref{l-bb} implies that $R$ is spanned by a movable rational curve.  Thus, by \cite[Theorem 4.10]{Kol91}, we  get the contraction as described in Case (2).

If $K_X+H$ is big, then we proceed with a step of a generalized minimal model program, given by a $K_X$-negative map $X\dasharrow X_1$ as described in \cite{HX13}, and we can replace $X$ by $X_1$. It follows from  termination of generalized flips (see \cite[Section 5]{HX13}) that the above process must terminate with one of the two cases of Theorem \ref{t-3mmp}.  
\end{proof}
\end{say}

\subsection{On the base point free theorem}\label{s-bpf}

The main aim of this section is to prove  Theorem \ref{t-abundance}. To this end, our main tool is the following result:

%
%
%
%
%

\begin{prop}\label{l-dim2} 
If $(X,B)$ is a log canonical threefold, where $K_X+B$ is nef and ${\rm Char}~ k=p>\frac{2}{a}$ where $a$ is the minimal nonzero coefficient of $B $. Assume $X$ has a dense open set $U$ which admits a dominant proper morphism $U\to V $ where $\dim (V)=2$. Assume also that $K_X+B$ is numerically trivial over the generic point of $V$.
Then $K_X+B$ is semiample. 
\end{prop}
\begin{proof}

 By the existence of resolution of singularities for curves and surfaces, we may assume that $\varphi$ induces a rational map $X\dashrightarrow Z$ where $Z$ is a smooth projective variety of dimension $n(L)$.  
 Let $\psi:Y\to X$ be a birational morphism which resolves the singularities of $X\dashrightarrow Z$, and whose existence is guaranteed by the main results in  \cites{Ab98,Cut04,CP08,CP09}. Thus, if we write $\psi^*(K_X+B)=K_{Y}+B_Y$, then $(Y,B_Y)$ is a sub log canonical pair.
 Let $f\colon Y\to Z$ be the induced map. 
  We can assume 
 \begin{enumerate} 
 \item  $f_*(\mathcal{O}_Y)=\mathcal{O}_Z$, and
\item  $f$ factors through an equidimensional morphism $Y^*\to Z$ (see \cite[Theorem 5.2.2]{RG71}), where $Y^*$ yields a morphism to $X$.
\end{enumerate}

\medskip 
We begin with the following Lemma.

\begin{lem}If $C$ is a normal complete curve defined over a field $\eta$ such that $\omega_C$ is anti-ample and $H^0(C,\mathcal{O}_C)=\eta$, then $C_{\bar{\eta}}$ is a conic in $\mathbb P^2_{\bar{\eta}}$. In particular, if ${\rm char}~ \eta >2$, then $C_{\bar{\eta}}\cong \mathbb{P}^1_{\bar{\eta}}$.

\end{lem}
\begin{proof}We have $H^0(C, \omega_{C})=0$ as $\omega^{-1}_{C}$ is ample.  
So the arithmetic genus of $C$ and $C_{\bar{\eta}}$ satisfy
$$a(C)=a(C_{\bar{\eta}})=0.$$
We know that $C_{\bar{\eta}}$ is irreducible. Let $C_{\bar{\eta}}^{\red}\subseteq C_{\bar{\eta}}$ be the reduced part and let $I$  be its ideal sheaf.  Then 
$$a(C^{\red}_{\bar{\eta}})\le a(C_{\bar{\eta}})=0,$$
which implies that $C^{\red}_{\bar{\eta}}$ is a smooth rational curve. 
Since for $j=0,1$, we have
$$H^j(C_{\bar{\eta}},\mathcal{O}_{C_{\bar{\eta}}})=H^j(C^{\red}_{\bar{\eta}},\mathcal{O}_{C^{\red}_{\bar{\eta}}})=H^j(C',\mathcal{O}_{C'})$$
for any $C^{\red}_{\bar{\eta}}\subseteq C'\subseteq C_{\bar{\eta}},$ we conclude that $H^j(C^{\red}_{\bar{\eta}},I^i/I^{i+1})=0,$ and in particular $I^i/I^{i+1}=\mathcal{O}_{\mathbb{P}^1}(-1)$ for any $i\le n$, where $n$ is the maximal non-negative integer such that $I^n\neq 0$. But then 
$$\omega_{C_{\bar{\eta}}}|_{C^{\red}_{\bar{\eta}}}\cong \mathcal{O}_{\mathbb{P}^1}(n-2).$$
Thus $n<2$, which implies that $C_{\bar{\eta}}$ is a conic in $\mathbb P^2_{\bar{\eta}}$, i.e., $C_{\bar{\eta}}$ is either a smooth rational curve or a planar double line and the latter case  can happen only if ${\rm char}~ \eta=2$. 
\end{proof}

The Lemma implies that if $\eta\in V$ is the general point then $X_{\bar{\eta}}\cong 
\mathbb{P}^1_{\bar{\eta}}$. Note that if $\eta$ is the general point of $Z$, then $Y_\eta$ is isomorphic to 
$X_\eta$. Denote by $B_{Y_\eta}$ the restriction $B_Y|_{Y_\eta}$.  Since by assumption $p>\frac{2}{a}$, 
where $a$ is the minimal non-zero coefficient of $B$ and since $K_{Y_\eta}+B_{Y_\eta}\sim_{\mathbb Q}0$, it 
follows that if $E$ is a  horizontal components of ${\rm Supp} ~B_Y$, then 
$$E\cdot Y_\eta\le \frac 1 a B\cdot Y_{\eta}=\frac 2 a<p.$$
In particular, $(Y_{\bar{\eta}},B_{Y_{\bar{\eta}}})$ is log canonical.
Moreover, if $E\to E' \stackrel {g_E}\longrightarrow Z$ is the Stein  factorisation of the morphism $E\to 
Z$, it follows that $\deg g_{E}<p$.

\medskip

We now define the $\mathbb Q$-divisor $D_{b} $ on $Z$ as the boundary part of  Kawamata subadjunction formula for $(Y,B_{Y})$ over $Z$. More precisely, for any prime divisor $W$ of $Z$, we define
 $$c_W=\sup \{t\mid (Y,B_Y+tf^*W ) \mbox{ is log canonical over the generic point of $W$}\}. $$
 Then $D_{b}:=\sum_W (1-c_W)W$ is a $\mathbb Q$-divisor on $Z$. 
 After possibly taking a log resolution, we may assume that $ (Z,{\rm Supp} D_b)$ is simple normal crossing \cite[Remark 7.7]{PS09}. Write $D_b=D_b^+-D_b^-$ where $D^+_b$ and $D^-_b$ are effective divisors and do not have common components. We fix  a divisor $\Gamma \ge D_b$, such that $(Z,\Gamma)$ is log canonical and the support of $\Gamma-D_b$ is contained in the negative part ${\rm Supp}~ D^-_b$ of $D_b$.

\medskip 

We now follow closely the arguments in \cite[Section 8]{PS09}. We denote by $\overline{\mathcal M}_{0,n}$ the moduli space of $n$-pointed stable curves of genus $0$ and we consider the 
universal family $\mathcal U_{0,n}\to \overline{\mathcal M}_{0,n}$. The varieties $\overline{\mathcal M}_{0,n}$ and $\mathcal 
U_{0,n}$ are both smooth and projective. We refer to \cite{Keel92} for a construction and  some basic properties of these 
varieties. In particular, the morphism  $g_n\colon \mathcal U_{0,n}\to \overline{\mathcal M}_{0,n}$ factors through a smooth projective 
variety $\overline {\mathcal U}_{0,n}$ such that the induced morphism $\sigma\colon \mathcal U_{0,n} \to \overline{\mathcal U}_{0,n}$
is a sequence of blow-ups with smooth centres and $\overline g_n\colon \overline{\mathcal U}_{0,n}\to \overline{\mathcal M}_{0,n}$ is a $\mathbb P^1$-bundle over $\overline{\mathcal M}_{0,n}$.

By taking a base change, we can assume that there is a diagram
\begin{diagram}
\mathcal{U}_{0,n}&\lTo& Y^2&\rTo&Y^1&\rTo^{h} & Y\\
\dTo^{g_n}&&\dTo^{f^2}& &\dTo & &\dTo^f\\
\overline{\mathcal{M}}_{0,n}&\lTo&Z^2&\rTo^{\mu}&Z^1&\rTo^{g} & Z 
\end{diagram}
such that:
\begin{enumerate}
\item $g$ is the composition of the morphisms $g_{E}$ defined above, for any horizontal component $E$ of ${\rm Supp}~B_Y$. Note that $\deg g_{E}< p$.  Let $Y^1$ be the normalization of the main component of $Y\times_Z Z^1$. In particular, if we define $B_{Y^1}$ by $K_{Y^1}+B_{Y^1}=h^*(K_Y+B_Y)$, then the horizontal components $E_1,\dots E_n$ of ${\rm Supp}~ B_{Y^1}$ correspond to  rational sections of $Z^1$,
\item Over the generic point $\eta_1$ of $Z^1$, we have that  $(Y^1,{\rm Supp} ~h^{-1}(B))|_{\eta_1}\in \overline{\mathcal M}_{0,n}(\eta_1)$, which yields a map $Z^1\dashrightarrow \overline{\mathcal M}_{0,n}$,  and
\item $\mu\colon Z^2\to Z^1$ is a birational morphism from a smooth surface $Z^2$ which resolves the singularities of the map $Z^1\dashrightarrow\overline{\mathcal{M}}_{0,n}$ and such that there exists a morphism  $f^2\colon Y^2\to Z^2$   which yields the morphisms $Y^2\to Y^1$ and $Y^2\to \mathcal{U}_{0,n}\times_{\overline{\mathcal{M}}_{0,n}}Z^2$. 
\end{enumerate}
 
\medskip 
  
Denote by $\rho\colon Y^2\to X$  the induced map. From the construction above, we easily get the following Lemma.

\begin{lem} Under the same assumptions as above, we have:

\begin{enumerate} 
\item Define the $\mathbb Q$-divisor $B_{Y^2}$ on $Y^2$ by $\rho^*(K_X+B)=K_{Y^2}+B_{Y^2}$. Then $(Y^2,B_{Y^2})$ is sub log canonical.
\item Let $D^2_{b}$ be the boundary part of the Kawamata subadjunction formula for $(Y^2,B_{Y^2})$ over $Z^2$. Then 
$(g\circ \mu)^*(K_{Z}+D_b)=K_{Z^2}+D_b^2.$
Furthermore, $(Z,D_b)$ is sub log canonical if and only if $(Z^2,D^2_b)$ is sub log canonical.  
\end{enumerate}
\end{lem}
\begin{proof}Since $g\colon Z^1\to Z$ is the composition of maps of degree less than $p$, it follows that the morphism  $\mu\colon Z^2\to Z$ is tamely ramified. Thus, we can apply the same arguments as in \cite[Proposition 5.20]{KM98}.
\end{proof}

We have the following canonical bundle formula:

\begin{lem}\label{l-bertini}
  Under the same assumptions as above, $(Z^2,D_b^2)$ is sub log canonical and
   there is a semiample divisor $D^2_m$ 
such that 
  $$(f^2)^*(K_{Z^2}+D^2_b+D^2_{m})\sim_{\mathbb{Q}} K_{Y^2}+B_{Y^2}.$$
\end{lem}
\begin{proof} Note that \cite[Theorem 2]{Kawamata97}  (see also \cite[Theorem 8.5]{PS09} and \cite[Section 3]{KMM94}) holds over any algebraically closed field, without any change in the proof. In particular, if $\mathcal P_1,\dots,\mathcal P_n$ are the sections of   $g_n\colon \mathcal U_{0,n}\to \overline{\mathcal M}_{0,n}$ corresponding to the marked points, and $d_1,\dots,d_n$ are the coefficients of $B_{Y^1}$ along $E_1,\dots,E_n$ respectively, then 
$$K_{\overline {\mathcal U}_n}+\sigma_* \mathcal D =  {\overline g}_n^*(K_{\overline{\mathcal M}_{0,n}} + L),$$
where  $\mathcal D=\sum_{i=1}^n d_i \mathcal P_i$ and $L$ is a semiample $\mathbb Q$-divisor on $\overline{\mathcal M}_{0,n}$.

It follows from the proof of  \cite[Theorem 8.1 ]{PS09} that there is a birational morphism 
$$j:Y^2\to \tilde{Y}^2:=  \overline{\mathcal{U}}_{0,n}\times_{\overline{\mathcal{M}}_{0,n}}Z^2$$ 
such that $h:\tilde{Y}^2\to Z^2$ is a $\mathbb{P}^1$-bundle over $Z^2$. Furthermore,   if we denote by $D^2_m$ the  pull-back of $L$ on $Z^2$, then 
$$h^*(K_{Z^2}+D^2_m)=K_{\tilde{Y}^2}+B^h_{\tilde{Y}^2},$$
where $B_{\tilde{Y}^2}=h_*(B_{Y^2})$ and $B^h_{{Y}^2}$ is the horizontal part of $B_{{Y}^2}$ over $Z^2$.

We claim that $j^*(K_{\tilde{Y}^2}+B_{\tilde{Y}^2})=K_{{Y}^2}+B_{{Y}^2}$. Assuming the claim,  $D^2_b$ can be computed on $(\tilde{Y}^2,B_{\tilde{Y}^2})$ by 
$$h^*D^2_b=B^v_{\tilde{Y}^2}=B_{\tilde{Y}^2}-B^h_{\tilde{Y}^2},$$ where $B^v_{\tilde{Y}^2}$ denotes the vertical part of $B_{{Y}^2}$ over $Z^2$.
 Thus, the Lemma easily follows. 

We now proceed with the proof of the claim. By the negativity lemma (cf. \cite[Lemma 1.17]{Kol13}), we have $$j^*(K_{\tilde{Y}^2}+B_{\tilde{Y}^2})-K_{{Y}^2}-B_{{Y}^2}=E\ge 0.$$ 
Since $K_{\tilde{Y}^2}+B_{\tilde{Y}^2}$ is the pull-back of a $\mathbb Q$-divisor on $Z$, we know that $-E$ is also nef over $Z$. But ${\rm Supp} ~E$ is exceptional over $Z$, i.e., for any codimension 1 point $P$ on $Z$ contained in $f({\rm Supp} ~E)$, we have that ${\rm Supp} E$ does not contain $f^{-1}(P)$. This implies that $E=0$. Thus, the claim follows. 
\end{proof}

\begin{rem} In general, in terms of the singularities, the canonical bundle formula in positive characteristic does not behave as well as in characteristic zero even for elliptic fibrations, due to the existence of fibrations with wild fibres (e.g. see \cite{BM77}). 
\end{rem}

We can pick an effective $\mathbb Q$-divisor $H^2\sim_{\mathbb{Q}} D^2_m$ and define $D_m$ to be $\frac{1}{\deg (g)}(\mu_*g_*H^2)$, such that  ${\supp D_m}$ does not have common components with ${\rm Supp(\Gamma)}$. In particular,
$(Z,\Gamma+D_m)$ is a log pair and $K_Y+B_Y\sim_{\mathbb{Q}}f^*(K_Z+D_b+D_m)$.

\begin{lem}\label{l-mini}
Under the same assumptions as above, if $m$ is a sufficiently divisible positive integer, then 
$$H^0(Z, m(K_Z+\Gamma+D_m))=H^0(Z,m(K_Z+D_b+D_m)).$$ 
\end{lem}

\begin{proof} We may write $B_Y=B^+_Y-B^-_Y+B^*_Y $, where $B^*_Y$  is the part of $B$ consisting of exactly  all the components of $B_Y$ which are exceptional over $Y^*$,  and the $\mathbb Q$-divisors $B^+_Y$ and $B^-_Y$ are effective and do not have common components.

 Let $G$ be the sum of all the prime divisors on $Y$ which are exceptional over $Y^*$ and let 
$$B'_Y=B^+_Y+{\rm Supp} ~B^-_Y+ tG$$ for some $t\gg 0$. Then since $B^-_Y$ and $B^*_Y$ are both exceptional over $X$, we have that for any sufficiently divisible positive integer $m$
$$H^0(Y,m(K_Y+B_Y))=H^0(X, m(K_X+B))=H^0(Y,m(K_Y+B'_Y)).$$
By the definition of the boundary part, for any sufficiently large $t$, we also have
$$f^*(K_Z+\Gamma+D_m)\le K_{Y}+B'_Y.$$
We conclude that, for sufficiently divisible positive integer $m$
$$\begin{aligned}
H^0(Y,m(K_Y+B'_Y))&\supseteq H^0(Z,m(K_Z+\Gamma+D_m))\\
&\supseteq H^0(Z, m(K_Z+D_b+D_m))\\
&=H^0(Y,m(K_Y+B_Y)).
\end{aligned}
$$
Thus, the claim follows. 
 \end{proof}

Since $Z$ is smooth and the coefficients of any prime divisor in $\Gamma+D_m$ is less than 1, by \cite{FT12}, we can run a MMP for $(Z,\Gamma+D_m)$,  which ends with a good minimal model $\pi: Z\to Z^m$. Let $F=\pi_*(\Gamma+D_m)$. 
It follows from Lemma \ref{l-mini} and the fact that $K_Z+D_b+D_m$ is nef that 
$$\pi^*(K_{Z^m}+F)=K_Z+D_b+D_m.$$ This implies that
$K_Y+B_Y\sim_{\mathbb{Q}}f^*(K_Z+D_b+D_m)$
is semiample. Thus we conclude that $K_X+B$ is semiample. 
\end{proof}

\begin{proof}[Proof of Theorem \ref{t-abundance}]
We first prove (1) and (2). 
Since it is sufficient to prove the statement of the Theorem after any base change of 
the ground field, we may assume that the ground field is uncountable.
By Theorem \ref{t-main2},  we can assume that $n(K_X+B)\le 2$. Thus, we only need to prove that if $n(K_X+B)=1$ or $2$, then $K_X+B$ is semiample. 

 Let $\varphi\colon U\to V$ be the nef reduction morphism of $K_X+B$ defined in Theorem \ref{t_nrm} where $U$ is an open subset of $X$. We distinguish two cases: 

\medskip 

\noindent {\bf Case 1: $n(K_X+B)=2$}. 

This case follows directly from Proposition \ref{l-dim2} and the fact that the restriction $(K_X+B)|_{U_{\eta}}$ on the generic fibre of $\varphi$ is numerically trivial.

\bigskip
\noindent{\bf Case 2: $n(K_X+B)=1$}.

We may assume that $\varphi$ induces a rational map $X\dashrightarrow Z$ where $Z$ is a smooth curve. Since $X$ is normal and $\dim Z=1$, the map $X\dashrightarrow Z$ is in fact a morphism which we denote by $h\colon X\to Z$.  
It suffices to show that $K_X+B\sim_{\mathbb{Q}}h^*G$ for some $\mathbb Q$-divisor $G$ on $Z$. Indeed, by Theorem \ref{t_nrm}, we have that $\deg G>0$ and the Theorem follows.

 Consider the Albanese morphism $a_X\colon X\to {\rm Alb}_X$, and denote by  $\phi\colon X\to S$  the Stein factorization of  the  morphism $(h,a_X)\colon X\to Z\times {\rm Alb}_X$. Denote by  $i\colon S \to Z\times {\rm Alb}_X$ the induced morphism. 
 Since the fibres of $Z$ are covered by rational curves, which are also mapped to points in ${\rm Alb}_X$, we know that $\dim (S)\le 2$. If $\dim(S)=1$, then there is an isomorphism $\rho\colon Z\to S$ such that $\rho\circ  h = \phi$.

We claim there exists a $\mathbb Q$-divisor  $H$ on $Z$ such that $K_X+B$ is numerically equivalent to $h^*H$. In fact, by the construction of the nef reduction map, we know that $(K_X+B)|_{K(Z)}$ is numerically trivial, where $K(Z)$ is the generic point of $Z$. Thus, there exists a $\mathbb Q$-divisor  $H$ on $Z$ and an effective $\mathbb Q$-divisor $E$ on $X$ such that 
$$K_X+B-h^*H\equiv E$$
where the support of $E$ is contained in a union of fibres of $h$ but  ${\rm Supp}(E)$ does not contain any fibre. Since $K_X+B$ is nef, it follows from Zariski's  lemma  that $E=0$, as claimed.

 Thus, if $n$ is a sufficiently large integer, we have  
 $$n(K_X+B-h^*H)\in ({\Pic}^0_{X/k})_{\red}=\Pic^0({\rm Alb}_X)$$ (see \cite[Remark 9.5.25 and Theorem 9.6.3]{Klei05}). Thus, we can find a divisor $M$ on ${\rm Alb}_X$ such that if $\pi\colon Z\times {\rm Alb}_X\to {\rm Alb}_X$ denotes the projection and $M'= (\pi\circ i)^*M$ then  $\deg M'=0$ and
 $$n(K_X+B-h^*H)\sim_{\mathbb{Q}} h^*(\rho^* M'),$$
 which implies $K_X+B\sim_{\mathbb{Q}}h^*(H+\frac{1}{n}\rho^* M'),$ as claimed. 
 
If $\dim S=2$,  then it follows directly from Proposition \ref{l-dim2} that $K_X+B$ is semiample. 
 
 \medskip 
We now proceed with the proof of (3).  
If $\kappa(X,K_X+B)=3$, then the result follows from 
\cite[Theorem 0.5]{Keel99}. Thus, by (1) and (2), it is enough to consider the case 
$n(X,K_X+B)=0$, which implies that $K_X+B$ is numerically trivial. Therefore, $K_X+B$ is $\mathbb{Q}$-linear 
equivalent to $0$ as this is true for any numerically trivial $\mathbb{Q}$-Cartier divisor on an algebraic 
variety defined over $\overline{\mathbb{F}}_p$.
\end{proof}


\bigskip

\end{document}